\documentclass[aap,MSNbibl,seceqn,dvips]{arximspdf}

\doi{10.1214/13-AAP933} 
\volume{24}
\issue{2}
\pubyear{2014}
\firstpage{679}
\lastpage{720}

\makeatletter
\newcommand{\rrvert}{\vert}
\newcommand{\llvert}{\vert}
\newtheorem{Theorem}{Theorem}[section]
\newproclaim{Definition}{Definition}[section]
\newtheorem{Proposition}{Proposition}[section]
\newtheorem{Lemma}{Lemma}[section]
\newtheorem{Corollary}{Corollary}[section]
\newproclaim{Remark}{Remark}[section]
\newproclaim{Example}{Example}[section]
\makeatother

\begin{document}
\begin{frontmatter}

\title{Runge--Kutta schemes for backward stochastic differential equations}
\runtitle{Runge--Kutta schemes for BSDEs}

\begin{aug}
\author[a]{\fnms{Jean-Fran\c{c}ois} \snm{Chassagneux}\corref{}\ead[label=e1]{j.chassagneux@imperial.ac.uk}\thanksref{t1}}
\and
\author[a]{\fnms{Dan} \snm{Crisan}\ead[label=e2]{d.crisan@imperial.ac.uk}\thanksref{t2}}
\runauthor{J.-F. Chassagneux and D. Crisan}
\affiliation{Imperial College London}
\address[a]{Department of Mathematics\\
Imperial College London\\
South Kensington Campus\\
London SW7 2AZ\\
United Kingdom\\
\printead{e1}\\
\phantom{E-mail: }\printead*{e2}} 
\end{aug}
\thankstext{t1}{Supported in part by Grant ANR-11-JS01-0007---LIQUIRISK.}
\thankstext{t2}{Supported in part by EPSRC Grant EP/H0005500/1.}

\received{\smonth{4} \syear{2012}}
\revised{\smonth{4} \syear{2013}}

%
\begin{abstract}
We study the convergence of a class of Runge--Kutta type schemes for
backward stochastic differential equations (BSDEs) in a Markovian
framework. The schemes belonging to the class under consideration
benefit from a certain stability property. As a consequence, the
overall rate of the convergence of these schemes is controlled by their
local truncation error. The schemes are categorized by the number of
intermediate stages implemented between consecutive partition time
instances. We show that the order of the schemes matches the number
$p$ of intermediate stages for $p\le3$. Moreover, we show that the
so-called order barrier occurs at $p=3$, that is, that it is not
possible to construct schemes of order $p$ with $p$ stages, when $p>3$.
The analysis is done under sufficient regularity on the final condition
and on the coefficients of the BSDE.
\end{abstract}

%

\end{frontmatter}


\section{Introduction}\label{sec1}

Let $(\Omega,{\mathcal F}, ({\mathcal F}_t)_{t \geq0},{\mathbb P})$ be
a filtered probability space endowed with an $({\mathcal F}_t)_{t \geq
0}$-adapted Brownian motion $(W_t)_{t \geq0}$. On $(\Omega,{\mathcal
F},\break  ({\mathcal F}_t)_{t \geq0},{\mathbb P})$ we consider the triplet $
( X,Y,Z ) = \{ ( X_{t},Y_{t},Z_{t} ),t\in[ 0,T ] \} $ of\break  $({\mathcal
F}_t)_{t \geq0}$-adapted stochastic processes satisfying the following
equations:
%
%
\begin{eqnarray}
X_{t} & =& X_{0} + \int_{0}^{t}b(X_{s})
\,\mathrm{d}s + \int_{0}^{t}\sigma
(X_{s})\,\mathrm{d} W_{s}, \label{eq sde X}
\\
Y_{t} &=& g(X_{T})+ \int_{t}^{T}
f(Y_{t},Z_{t}) \,\mathrm{d}t -\int_{t}^{T}
Z_{t} \,\mathrm{d}W_{t}. \label{eq bsde YZ}
\end{eqnarray}
System (\ref{eq sde X})--(\ref{eq bsde YZ}) is called a (decoupled)
forward-backward stochastic differential equation (FBSDE).

The process $X$, called the forward component of the FBSDE, is a
$d$-dimensional diffusion satisfying a stochastic differential equation
(SDE) with Lipschitz-continuous coefficients $b\dvtx
\mathbb{R}^{d}\rightarrow\mathbb{R}^{d}$ and $\sigma\dvtx
\mathbb{R}^{d}\rightarrow\mathbb{R}^{d}\times\mathbb{R}^{d}$.

The pair of processes $(Y, Z)$ satisfy the backward
stochastic differential equation (BSDE) (\ref{eq bsde YZ}).
The process $Y$ is a one-dimensional stochastic
process with final condition $Y_{T}=g (X_{T})$, where $g\dvtx\mathbb
{R}%
^{d}\rightarrow\mathbb{R}$ is a differentiable function with
continuous and bounded first derivative [i.e., $g \in
C^1_b(\mathbb{R}^{d})$].
The process $Z=(Z^1,\ldots,Z^d)$ is a $d$-dimensional process,
written, by
convention, as a row vector.
The function $%
f\dvtx \mathbb{R}\times\mathbb{R}%
^{d}\rightarrow\mathbb{R}$ referred to as ``the driver,'' is assumed to
be Lipschitz
continuous.\tsup{3,4}\setcounter{footnote}{3}\footnotetext{These
assumptions will be strengthened in the following
section.}\setcounter{footnote}{4}\footnotetext{For the reader's
convenience, we only consider drivers depending on $Y$ and $Z$;
however, the results and the analysis provided here apply to drivers
depending also on $X$.}

The existence and uniqueness of solutions of system (\ref{eq sde
X})--(\ref{eq bsde YZ}) was first addressed by Pardoux and Peng in
\cite{parpen90}. Since then, a large number of papers have been
dedicated to the study of FBSDEs. In particular, it is well known that
under the Lipschitz-continuity assumption of the coefficients, the
following estimate holds true:
%
%
\begin{equation}
\label{eq basic estimates} \quad\mathbb{E} \Bigl[\sup_{t \in[0,T]}|X_t|^p
\Bigr] + \mathbb{E} \biggl[\sup_{t
\in[0,T]}|Y_t|^2
+ \int_{0}^T|Z_s|^2 \,
\mathrm{d}s \biggr] < \infty\qquad\forall p > 0.
\end{equation}

Moreover, Pardoux and Peng showed in \cite{parpen92} that
\[
Y_t=u(t,X_t),\qquad Z_t=\nabla
u^{ \top} (t,X_t)\sigma(X_t),\qquad t\in[0,T],
\]
where $u\in C^{1,2}([0,T]\times\mathbb{R}^d)$ is the solution of the
final value Cauchy problem
%
%
\begin{eqnarray}
\quad L^{(0)}u(t,x)&=&-f \bigl( u(t,x),\nabla u^{ \top} (t,x)
\sigma( x ) \bigr),\qquad t\in\lbrack0,T), x\in\mathbb{R}^{d},
\label{b_pde_t}
\\
u(T,x)&=&g (x),\qquad x\in\mathbb{R}^{d} \label{b_pde_T}
\end{eqnarray}
with $L^{(0)}$ defined to be the second order differential
operator
%
%
\begin{equation}
\label{operator} L^{(0)}=\partial_{t} + \sum
_{i=1}^{d}b_i \partial_{x_i}+
\frac12 \sum_{i=1}^{d}a_{ij}
\partial_{x_i}\partial_{x_j}
\end{equation}
and $a=(a_{ij})=\sigma\sigma^{\top}$.

There is a vast literature dedicated to the approximation of solutions
to stochastic differential equations. In particular, obtaining
approximations of the distribution of the forward component $X$ has
been largely resolved in the last thirty years. One can refer to
\cite{klopla92} and the references therein for a systematic study of
numerical methods for approximating $X$. Such methods are classical by
now. More recently, Kusuoka, Lyons, Ninomiya and Victoir
\cite{kus01,lyovic04,nim03a,nim03b,nimvicXX} developed several
numerical algorithms for approximating $X$ based on Chen's iterated
integrals expansion. These new algorithms generate\vadjust{\goodbreak} an approximation
of
the solution of the SDE in the form of the empirical distribution of a
cloud of particles with deterministic trajectories.

By comparison, there are very few numerical methods for approximating
the backward component. In this paper, we introduce a large class of
numerical schemes for approximating solutions of BSDEs. These schemes
are based on the well-known Runge--Kutta methods for ODEs and include
new high order schemes as well as existing low order schemes such as
the classical extension of the Euler scheme to BSDEs; see, for example,
\cite{boutou04,boueli08,criman10,goblab07}.

The approximations presented below are associated to an
arbitrary, but fixed, partition $\pi$ of the interval $[0,T]$, $\pi=
\{t_{0} = 0 <\cdots<t_i<t_{i+1}<\cdots< t_{n} = T\}$. We denote $h_{i}=
t_{i +1}-t_{i }$, $i=0,\ldots,n-1$ and $|\pi| = \max_{i}h_{i}$.
Let $(Y_i,Z_i)$ be the approximation of $(Y_{t_i},Z_{t_i})$ for
$i=1,\ldots,n$. The construction of the approximating process is done
in a recursive manner, backwards in time. We describe in the following
the salient features of the class of approximations considered in this
paper.


%
%
\begin{Definition}\label{de RK scheme}
\begin{longlist}[(ii)]
\item[(i)] The terminal condition is given by the pair
$(Y_n,Z_n)=(g(X_T),\break  \nabla g^\top(X_T)\*\sigma(X_T))$.

\item[(ii)] For $i\le n-1$, the transition from $(Y_{i+1},Z_{i+1})$ to
$(Y_i,Z_i)$ involves $q$~stages, with $q\ge1$.
Given $q+1$ positive coefficients $0=: c_1 < c_2 \le\cdots\le c_j \le
\cdots c_{q} \le c_{q+1}:= 1$, we introduce the intermediate
``instances'' of computation $t_{i,j}:= t_{i +1}-c_j h_{i}, $ and
define $(Y_{i,j},Z_{i,j})$, $j=1,\ldots,q+1$ as follows: by convention,
$(Y_{i,1},Z_{i,1}) = (Y_{i +1},Z_{i +1})$ and $(Y_{i,q+1},Z_{i,q+1} ) =
(Y_{i },Z_{i })$. Then, for $1<j\le q$,
%
%
\begin{eqnarray}
\label{eq RK scheme general} Y_{i,j} &=& \mathbb{E}_{t_{i,j}}
\Biggl[ Y_{i +1} + c_jh_{i}\sum
_{k=1}^{j} a_{jk} f(Y_{i,k},Z_{i,k})
\Biggr],
\\
Z_{i,j} &=& \mathbb{E}_{t_{i,j}} \Biggl[H^i_{j}
Y_{i +1} + h_{i}\sum_{k=1}^{j-1}
\alpha_{jk} H^i_{j,k} f(Y_{i,k},Z_{i,k})
\Biggr].
\end{eqnarray}

Finally, the approximation at step~($\mathrm{i}$) is given by
%
%
\begin{eqnarray}
\label{eq RK scheme stability last} Y_{i } &=& \mathbb{E}_{t_{i }}
\Biggl[ Y_{i +1} + h_{i}\sum_{j=1}^{q+1}
b_{j}f(Y_{i,j},Z_{i,j}) \Biggr],
\\
Z_{i } &=& \mathbb{E}_{t_{i }} \Biggl[H^i_{q+1}
Y_{i +1} + h_{i}\sum_{j=1}^{q}
\beta_{j} H^i_{q+1,j} f(Y_{i,j},Z_{i,j})
\Biggr].
\end{eqnarray}
%

The coefficients $(a_{jk})_{1\le j,k \le q}$,
$(\alpha_{jk})_{1\le j,k \le q}$, $(b_{j})_{1\le j \le q+1}$ and
$(\beta_{j})_{1\le j \le q}$ take their values in $\mathbb{R}$ with $a_{1j}$,
$\alpha_{1j}$, $1 \le j \le q$ and $a_{jk}$, $\alpha_{jk}$, $1\le j<k
\le q$ set to $0$. Moreover, the following holds:
%
%
\begin{equation}
\label{eq coef order 0} \sum_{k=1}^j
a_{jk} = \sum_{k=1}^{j-1}
\alpha_{jk}\mathbf{1}_{\{ c_k<c_j\}} = c_j,\qquad j \le q.
\end{equation}
The random variables $H^i_j$, $H^i_{j,k}$, $k \le j$ are
$\mathcal{F}_{t_{i,j}}$-measurable, for all $j\le q+1$, $i < n$ and
have the property that, for all $1\le k < j \le q+1$, $i< n$,
%
%
\begin{eqnarray}
\label{eq prop H} \qquad\mathbb{E}_{t_{i,j}} \bigl[H^i_j
\bigr]&=& \mathbb{E}_{t_{i,j}} \bigl[H^i_{j,k} \bigr]
= 0\quad\mbox{and}\quad\mathbb{E} \bigl[h_{i} \bigl|H^i_j
\bigr|^2 \bigr]+ \mathbb{E} \bigl[h_{i} \bigl|H^i_{j,k}
\bigr|^2 \bigr] \le\Lambda,
\end{eqnarray}
where $\Lambda$ is a positive constant which does not depends on $\pi$.
\end{longlist}
\end{Definition}

Observe that $Y_n$, $Z_n$ belong to $ L^2(\mathcal{F}_{t_n})$, where
for $t \in[0,T]$, $L^2(\mathcal{F}_t)$ is the space of $\mathcal
{F}_t$-measurable random variables $U$ such that $\mathbb{E} [|U|^2
]<\infty$. This is an immediate consequence of estimates (\ref{eq basic
estimates}) and the fact that $g \in C^1_b$. Moreover, an easy
(backward) induction proves that the schemes are well defined for
$|\pi|$ small enough and that $Y_i$, $Z_i$ belong to
$L^2(\mathcal{F}_{t_{i }})$ for all $i \le n$.

In the sequel, we will refer to the schemes defined above by specifying
the \mbox{$H$-}coefficients and using the following tableau for the
other coefficients:
\[
\begin{array} {c@{\hspace*{5pt}}|@{\hspace*{5pt}}c@{\quad}c@{\quad}c@{\quad}c@{\hspace*{5pt}}|@{\hspace*{5pt}}c@{\quad}c@{\quad}c} c_1 = 0 & a_{11} &
\cdots& a_{1q} & 0 & \alpha_{11}& \cdots&
\alpha_{1q}
\\
\vdots& \vdots& &\vdots& \vdots& \vdots& &\vdots
\\
c_j& a_{j1} & \cdots& a_{jq} & 0 &
\alpha_{j1} &\cdots& \alpha_{jq}
\\
\vdots& \vdots& & \vdots& \vdots& \vdots& &\vdots
\\
c_q& a_{q1} & \cdots& a_{qq} & 0 &
\alpha_{q1}& \cdots& \alpha_{qq}
\\
&&&&
\\[-12pt]
\hline
c_{q+1} = 1 & b_{1} & \cdots&
b_{q} & b_{q+1} & \beta_{1} &\cdots&
\beta_{q} \end{array}.
\]


This notation is a natural extension of the classical notation used in
the ODEs framework; see, for example, \cite{but08}.

If the scheme is explicit for the last stage, that is, $b_{q+1}=0$, we
will omit this column in the coefficients tableau. We will also
generally omit the ``0'' coefficients in the tableau and use ``*'' to
denote a coefficient whose value is arbitrary.

Finally, let us also introduce for later use
%
%
\begin{equation}
\label{eq de tilde alpha-beta} \tilde{\alpha}_{jk} =
\alpha_{jk} \mathbf{1}_{\{c_k <
c_j\}}\quad\mbox{and}\quad\tilde{
\beta_j} = \beta_j \mathbf{1}_{\{c_j
< 1\}}.
\end{equation}

\subsection{General formulation of one-step schemes}\label{subse ge res}
It is convenient to rewrite the approximations
defined above in a more general setting as follows.

%
%
\begin{Definition}[(One-step scheme)]\label{de one-step
scheme}
\begin{longlist}[(ii)]
\item[(i)] The terminal condition is given by a pair $(Y_n,Z_n) \in
L^2(\mathcal{F}_{T})$.

\item[(ii)] For $i\le n-1$, the transition from $(Y_{i+1},Z_{i+1})$ to
$(Y_i,Z_i)$ is given by
%
%
\begin{equation}
\label{eq one-step scheme} \cases{ Y_{i } = \mathbb{E}_{t_{i }}
\bigl[ Y_{i +1} + h_{i}\Phi_i^Y(t_{i +1},Y_{i
+1},Z_{i +1},h_{i})
\bigr], \vspace*{3pt}
\cr
Z_{i } = \mathbb{E}_{t_{i }}
\bigl[H^{i}_{q+1}Y_{i +1} + h_{i}\Phi
_i^Z(t_{i +1},Y_{i +1},Z_{i
+1},h_{i})
\bigr],}\vadjust{\goodbreak}
\end{equation}
where $\Phi_i^Y$, $\Phi_i^Z$ are functions from $\mathbb{R}_+ \times
L^2(\mathcal{F}_{{t_{i+1}} }) \times L^2(\mathcal{F}_{{t_{i+1}} })
\times\mathbb{R}_+^* $ to $L^2(\mathcal{F}_{t_{i +1}})$, $0 \le i \le
n-1$.
\end{longlist}
\end{Definition}

%
%
\begin{Remark}
In the case of the scheme given in Definition~\ref{de RK scheme}, the
functions $\Phi_i^Y, \Phi_i^Z$ depend implicitly of the coefficients
$(a_{jk})_{1\le j,k \le q}$, $(\alpha_{jk})_{1\le j,k \le q}$,
$(b_{\hspace*{-0.5pt}j\hspace*{-0.2pt}})_{\hspace*{-0.1pt}1\le j \le
q+1\hspace*{-0.2pt}}$ and
$(\beta_{\hspace*{-0.4pt}j\hspace*{-0.2pt}})_{\hspace*{-0.1pt}1\le j
\le q\hspace*{-0.2pt}}$ and the random variables
$(H^i_{\hspace*{-0.2pt}j})_{\hspace*{-0.1pt}1\le j \le
q+1\hspace*{-0.2pt}}$,
$(H^i_{\hspace*{-0.2pt}j,k})_{\hspace*{-0.1pt}1\le j,k \le q\hspace*{-0.2pt}}$. 
\end{Remark}

\subsubsection{Order of convergence}\label{sec1.1.1}
The global error we investigate here is given by the pair
$(\mathcal{E}_Y(\pi),\mathcal{E}_Z(\pi))$, where
%
\begin{eqnarray*}
\mathcal{E}_Y(\pi)&:=& \max_{0\le i \le n} \mathbb{E}
\bigl[|Y_{t_{i }}-Y_i|^2 \bigr],\qquad 
\\
\mathcal{E}_Z(\pi)&:=& \sum_{i=0}^{n-1}
h_{i}\mathbb{E} \bigl[|Z_{t_{i }}-Z_i|^2
\bigr].
\end{eqnarray*}
To control these errors we will use \textit{the local truncation error}
for the pair $(Y,Z)$ defined as
%
%
\begin{eqnarray}
\label{eq de trunc error loc YnZ} \qquad\eta_i &:=&
\eta^Y_i + \eta^Z_i,\qquad \bigl(
\eta^Y_i, \eta^Z_i \bigr):=
\biggl(\frac1{h_{i}^2} \mathbb{E} \bigl[|Y_{t_{i }}
- \hat{Y}_{t_{i }}|^2 \bigr],\mathbb{E} \bigl[|Z_{t_{i }}-
\hat{Z}_{t_{i }}|^2 \bigr] \biggr)
\end{eqnarray}
with
%
%
\begin{equation}
\label{eq de hY hZ} \cases{ \hat{Y}_{t_{i
}}:=\mathbb{E}_{t_{i }}
\bigl[Y_{{t_{i+1}} } + h_i \Phi_i^Y(t_{i
+1},Y_{t_{i +1}},Z_{t_{i +1}},h_{i})
\bigr], \vspace*{3pt}
\cr
\hat{Z}_{t_{i }}:= \mathbb{E}_{t_{i }}
\bigl[H^{i}_{q+1}Y_{{t_{i+1}} } + h_{i}
\Phi_i^Z(t_{i +1},Y_{t_{i +1}},Z_{t_{i +1}},h_{i})
\bigr].}
\end{equation}

The \textit{global truncation error} for a given grid $\pi$ is given by
%
%
\begin{eqnarray}
\label{eq de trunc error YnZ} \mathcal{T}(\pi)&:=& \mathcal{T}_Y(\pi
) + \mathcal{T}_Z( \pi),
\nonumber
\\[-4pt]
\\[-12pt]
\bigl(\mathcal{T}_Y(\pi), \mathcal{T}_Z(\pi) \bigr)&:=&
\Biggl( \sum_{i=0}^{n-1} h_{i}
\eta^Y_i, 
\sum
_{i=0}^{n-1} h_{i} \eta^Z_i
\Biggr),
\nonumber
\end{eqnarray}
where $\mathcal{T}_Y$ is the global truncation error for $Y$, and
$\mathcal{T}_Z$ is the global truncation error for $Z$ defined as
above.

The main results of the paper refer to the rate of convergence of the
various approximations belonging to the class described in
Definition~\ref{de RK scheme}.

%
%
\begin{Definition}\label{de order}
An approximation is said to have a \textit{global truncation error} of
order $m$ if we have
\[
\mathcal{T}(\pi) \le C |\pi|^{2m}
\]
for all sufficiently smooth\footnote{The required regularity
assumptions will be stated in the theorems below.} solutions to
(\ref{b_pde_t})--(\ref{b_pde_T}) and all partitions $\pi$ with
sufficiently small mesh size.
\end{Definition}

%
%
\begin{Remark}
Observe that we consider the sum of the global truncation error for the
$Y$ component and the $Z$ component to define the \textit{order} of an
approximation. It is clear that if one considers BSDEs where the driver
$f$ depends only on $Y$ and is only interested in the error on the $Y$
part, it would be more judicious to use only $\mathcal{T}_Y$ in the
definition of the order of the method. But our goal here is to deal
with the most general case, where $f$ depends on both $Y$ and $Z$.
\end{Remark}
%

\subsubsection{Stability}\label{sec1.1.2}
To connect the truncation error with the global approximation error, we
introduce the notion of $L^2$-stability for the schemes given in
Definition~\ref{de one-step scheme}. By \textit{stability} we
mean---roughly speaking---that the outcome of the scheme is
``reasonably'' modified if we ``reasonably'' perturb the scheme.


We thus introduce a perturbed scheme,
%
%
\begin{equation}
\label{eq RK explicit scheme} \cases{ \tilde{Y}_{i } =
\mathbb{E}_{t_{i }} \bigl[\tilde{Y}_{i +1} + h_{i}
\Phi^Y(t_{i },h_{i},\tilde{Y}_{i +1},
\tilde{Z}_{i +1}) + \zeta^{Y}_{i
} \bigr],
\vspace*{3pt}
\cr
\tilde{Z}_{i } = \mathbb{E}_{t_{i }}
\bigl[H^i_{q+1} \tilde{Y}_{i +1} + h_{i}
\Phi^Z( t_{i },h_{i},\tilde{Y}_{i +1},
\tilde{Z}_{i +1},h_{i}) + \zeta^{Z}_{i }
\bigr],}
\end{equation}
where $\zeta^{Y}_{i }$, $\zeta^{Z}_{i }$ belongs to $L^2(\mathcal
{F}_{{t_{i+1}} })$, for all $i<n$ and with terminal values
$\tilde{Y}_n$~and~$\tilde{Z}_n$ belonging to $L^2(\mathcal{F}_{T})$.

For $0 \le i \le n$, we denote $\delta Y_{i }:= Y_{i } - \tilde{Y}_{i }$
and $\delta Z_{i }:= Z_{i } - \tilde{Z}_{i }$ and consider the following
definition of stability.

%
%
\begin{Definition}[($L^2$-Stability)]\label{de L2-stability}
The scheme given in Definition~\ref{de one-step scheme} is said to be
$L^2$-stable if
\begin{eqnarray*}
&&\max_{i}\mathbb{E} \bigl[|\delta Y_{i }|^{2}
\bigr] + \sum_{i=0}^{n-1}h_i
\mathbb{E} \bigl[|\delta Z_{i }|^{2} \bigr]
\\
&&\qquad\le C \Biggl( \mathbb{E} \bigl[|\delta Y_{n}|^{2}
+h_{n-1} |\delta Z_{n}|^{2} \bigr] + \sum
_{i=0}^{n-1} h_{i}\mathbb{E}
\biggl[\frac{1}{h_{i}^{2}} \bigl|\mathbb{E}_{t_{i }} \bigl[\zeta^{Y}_{i }
\bigr] \bigr|^{2} + \bigl|\mathbb{E}_{t_{i }} \bigl[\zeta^{Z}_{i }
\bigr] \bigr|^{2} \biggr] 
\Biggr)
\end{eqnarray*}
for all sequences $\zeta^{Y}_{i }$, $\zeta^{Z}_{i }$ of $L^2(\mathcal
{F}_{{t_{i+1}} })$-random
variables and terminal values $(Y_n, Z_n)$, $(\tilde{Y}_n,
\tilde{Z}_n)$ belonging to $L^2(\mathcal{F}_T)$.
\end{Definition}

Under a reasonable assumption on the functions $\Phi_i^Y$ and
$\Phi_i^Z$, $i \le n-1$, introduced in (\ref{eq one-step scheme}), we
are able to prove the stability of the schemes given in
Definition~\ref{de one-step scheme}.

%
%
\begin{Theorem}[(Sufficient condition for $L^2$-stability)]\label{th stab}
Assume that, for some given grid $\pi$ and for $i\le n-1$, we have
%
%
%
\begin{eqnarray}
&&\mathbb{E}_{t_{i }} \bigl[ \bigl| \Phi_i^Y(
t_{i +1},U,V,h_{i}) - \Phi_i^Y(t_{i +1},
\tilde U, \tilde V,h_{i}) \bigr|^2 \bigr]
\nonumber
\\[-8pt]
\label{eq prop PhiY}
\\[-8pt]
&&\qquad\le C \biggl( \frac1h_{i} \bigl(\mathbb{E}_{t_{i }}
\bigl[|\delta U|^2 \bigr] - \bigl|\mathbb{E}_{t_{i }} [\delta U ]
\bigr|^2 \bigr) + \mathbb{E}_{t_{i }} \bigl[|\delta
U|^2 + |\delta V|^2 \bigr] \biggr),
\nonumber
\\
&&\mathbb{E}_{t_{i }} \bigl[ \bigl| \Phi_i^Z(
t_{i +1},U,V,h_{i}) - \Phi_i^Z(t_{i +1},
\tilde U, \tilde V,h_{i}) \bigr|^2 \bigr]
\nonumber
\\[-8pt]
\label{eq prop PhiZ}
\\[-8pt]
&&\qquad\le\frac{C}{h_{i}} \biggl( \frac1h_{i} \bigl(\mathbb
{E}_{t_{i }} \bigl[|\delta U|^2 \bigr] - \bigl|
\mathbb{E}_{t_{i }} [\delta U ] \bigr|^2 \bigr) +
\mathbb{E}_{t_{i }} \bigl[|\delta U|^2 +|\delta
V|^2 \bigr] \biggr),
\nonumber
\end{eqnarray}
where $U$, $V$, $\tilde{U}$, $\tilde{V}$ belong to
$L^2(\mathcal{F}_{t_{i +1}})$, $\delta U:= U - \tilde U$ and $\delta V:=
V - \tilde V$, then the scheme in Definition~\ref{de one-step scheme}
is $L^2$-stable.
\end{Theorem}

The following proposition connects the truncation error with the
approximation error.

%
%
\begin{Proposition}\label{pr pre-conv}
Assume that the functions $\Phi_i^Y$ and $\Phi_i^Z$ satisfy (\ref{eq
prop PhiY})--(\ref{eq prop PhiZ})
and $(Y_n,Z_n)=(g(X_T), \nabla g^\top(X_T)\sigma(X_T))$. Then there
exists a constant $C$~independent of the partition $\pi$ such that 
%
%
\begin{eqnarray}
\mathcal{E}_Y(\pi) + \mathcal{E}_Z(\pi) &\le& C
\mathcal{T}(\pi).
\end{eqnarray}
\end{Proposition}

The proofs of Theorem~\ref{th stab} and Proposition~\ref{pr pre-conv}
are postponed to the \hyperref[sec7]{Appendix}.

\subsubsection{Convergence results}\label{sec1.1.3}

As an application of Definitions~\ref{de order}~and~\ref{de
L2-stability}, and Proposition~\ref{pr pre-conv}, we state the
following general convergence results (the proofs are postponed to the
\hyperref[sec7]{Appendix}):

%
%
\begin{Proposition} \label{pr ge conv res}
If the method is of order $m$
and $\Phi_i^Y$ and $\Phi_i^Z$ satisfy \mbox{(\ref{eq prop PhiY})--(\ref
{eq prop PhiZ})} and $(Y_n,Z_n)=(g(X_T), \nabla
g^\top(X_T)\sigma(X_T))$, then there exists a constant $C$ independent
of the partition $\pi$ such that
%
%
\begin{equation}
\mathcal{E}_Y(\pi) + \mathcal{E}_Z(\pi) \le C|
\pi|^{2m}.
\end{equation}
\end{Proposition}

Let us conclude this section with the main case of interest for us
here, namely the Runge--Kutta schemes given in Definition~\ref{de RK
scheme}.

%
%
\begin{Theorem} \label{th conv RK scheme}
\textup{(i)}~For the schemes given in Definition~\ref{de RK scheme}, if
$f$ is Lipschitz-continuous, we have that the functions $\Phi_i^Y$ and
$\Phi_i^Z$ satisfy (\ref{eq prop PhiY})--(\ref{eq prop PhiZ}) provided
$|\pi|$ is small enough. As a result, the schemes are $L^2$-stable.

\textup{(ii)} Moreover, if the method is of order $m$, then we have
%
%
\begin{equation}
\mathcal{E}_Y(\pi) + \mathcal{E}_Z(\pi) \le C|
\pi|^{2m},
\end{equation}
provided $|\pi|$ is small enough.
\end{Theorem}

%
%
\begin{Remark}
In this paper, we are only interested in obtaining an upper bound for
the global approximation error $\mathcal{E}_Y(\pi) + \mathcal
{E}_Z(\pi)$, in terms of
$|\pi|$. An asymptotic expansion of this error in term of $|\pi|$ would
also be of interest as it may lead to the use of Romberg--Richardson's
extrapolation method. This work is left for future research.
\end{Remark}

\subsection{Order of convergence of Runge--Kutta methods}\label{sec1.2}

It is a nontrivial task to classify the approximations belonging to the
class described by Definition~\ref{de RK scheme} through their order of
convergence. The order of convergence of a particular scheme depends on
several factors. First, it will depend on the number of intermediate
steps it uses. Moreover, up to a certain level, the higher the
smoothness of the pair $(u,f)$, the better the order is. However, there
is a level of smoothness beyond which the order of approximation cannot
typically be improved. This level is identified below through the
condition $(\mathbf{H}r)_p$, where $p=1,2,\ldots$ is the number of
intermediate steps required by the approximation. We show below that,
provided the underlying framework satisfies a certain nondegeneracy
condition called $(\mathbf{H}o)_p$, the order of the approximation
cannot be improved through additional smoothness. This is achieved by
identifying the leading order term in the expansion of the error of the
approximation. However, should this leading order term be equal to
zero, the order of the approximation will be higher. The analysis of
the leading error term tells us that, for example, if the driver
satisfies the additional constraint $f^z=0$ (i.e., it is independent of
$Z$, $f^z$ denoting the partial derivative of $f$ with respect to $z$),
then there are two-stage schemes of order three. However, if $f^z
\ne0$, then two-stage schemes will typically have order
two.\looseness=-1

\subsubsection{Smoothness and nondegeneracy assumptions}\label{sub se ass}
We study the order of the methods given in Definition~\ref{de RK
scheme} using It\^o--Taylor expansions \cite{klopla92}. This requires
the smoothness of the value function $u$. In order to state precisely
these assumptions, we recall some notations of Chapter~5 (see
Section~5.4) in \cite{klopla92}.

Let
\[
\mathcal{M}:= \{\oslash\}\cup\bigcup_{m =1}^\infty
\{0,\ldots,d\}^m
\]
be the set of multi-indices with entries in $\{0,\ldots,d\}$ endowed
with the measure $\ell$ of the length of a multi-index
[$\ell(\oslash)=0$ by convention].

We introduce the concatenation operator $*$ on $\mathcal{M}$ for
multi-indices with finite length $\alpha= (\alpha_1,\ldots,\alpha_p)$,
$\beta= (\beta_1,\ldots,\beta_q)$ then $\alpha*\beta=
(\alpha_1,\ldots,\alpha_p,\break \beta_1,\ldots,\beta_q)$.

For a multi-index $\alpha$ with positive finite length, we write
$-\alpha$ (resp., $\alpha-$) the multi-index obtained by deleting the
first (resp., last) component of $\alpha$.
On the set~$\mathcal{M}$, let $n(\alpha)$ be the number of zero in a
multi-index $\alpha$ with finite length.\vadjust{\goodbreak}

Given a multi-index $\alpha$, we denote by $\alpha^+$ the multi-index
obtained from $\alpha$ by deleting all its zero components.

For $j \in\{0,1,\ldots,d\}$, we denote by $(j)_m$ the multi-index with
length $m$ and whose entries are all equal to $j$.

A nonempty subset $\mathcal{A}\subset\mathcal{M}$ is called a
hierarchical set if
\[
\sup_{\alpha} \ell(\alpha) < \infty\quad\mbox{and}\quad-\alpha
\in\mathcal{A}\qquad\forall\alpha\in\mathcal{A}\setminus\{ \oslash\}.
\]

For any hierarchical $\mathcal{A}$ set, we consider the remainder set
$\mathcal{B}(\mathcal{A})$ given by
\[
\mathcal{B}(\mathcal{A}):= \{\alpha\in\mathcal{M}\setminus \mathcal{A}| -\alpha
\in\mathcal{A}\}.
\]

We will use in the sequel the following sets of multi-indices, for
$n\ge0$:
\[
\mathcal{A}_{n}:= \bigl\{\alpha| \ell(\alpha) \le n \bigr\}
\]
and observe that $\mathcal{B}(\mathcal{A}_{n})=\mathcal
{A}_{n+1}\setminus\mathcal{A}_{n}$.

For $j \in\{1,\ldots,d\}$, we consider the operators
\[
L^{(j)} = \sum_{k=1}^d
\sigma^{kj} \partial_{x_k}.
\]

For a multi-index $\alpha=(\alpha_1,\ldots,\alpha_p)$, the
iteration of
these operators has to be understood in the following sense:
\[
L^\alpha:= L^{(\alpha_1)}\circ\cdots\circ L^{(\alpha_p)}.
\]
By convention, $L^{\oslash}$ is the identity operator; recall also the
definition of the operator $L^{(0)}$ given in (\ref{operator}). One
can observe that $L^{\alpha*\beta}=L^\alpha\circ L^\beta$.

Let $C^k_b$ be the set of all $k$-times continuously differentiable
functions with all partial derivatives bounded. For a multi-index with
finite length $\alpha$, we consider the set $\mathcal{G}^\alpha$ of all
functions $v\dvtx[0,T]\times\mathbb{R}^d \rightarrow\mathbb{R}$ for
which $L^\alpha v$ is well defined and continuous. We also introduce
$\mathcal{G}^\alpha_b$ the subset of all functions $v
\in\mathcal{G}^\alpha$ such that the function $L^\alpha v$ is bounded.
For $v \in\mathcal{G}^\alpha$, we denote $L^\alpha v$ by $v^\alpha$.

Finally, for $n\ge1$, we define the set $\mathcal{G}^n_b$ of function
$v$ such that $v \in\mathcal{G}^\alpha_b$ for all $\alpha\in\mathcal
{A}_n\setminus \{\oslash\}$.

We are now ready to state the smoothness assumption on the
value function $u$ we shall use:
\begin{longlist}[$(\mathbf Hr)_1$]
\item[$(\mathbf Hr)_1$] The value function $u$ belongs to
$\mathcal{G}^2_b$ and $f \in
C^1_b$.\vspace*{2pt}

\item[$(\mathbf Hr)_2$] The value function $u$ belongs to
$\mathcal{G}^3_b$ and $f \in C^2_b$.\vspace*{2pt}

\item[$(\mathbf Hr)_3$] The value function $u$ belongs to
$\mathcal{G}^4_b$ and $f \in C^3_b$.\vspace*{2pt}

\item[$(\mathbf Hr)_4$] The value function $u$ belongs to
$\mathcal{G}^5_b$ and $f \in C^5_b$.
\end{longlist}

Instead of making assumptions on the coefficient $b$ and
$\sigma$, we shall use in the sequel the following ``nondegeneracy''
assumption when stating the necessary order conditions:
\begin{longlist}[$(\mathbf Ho)_2$]
\item[$(\mathbf Ho)_1$] There exists some function $g
\in\mathcal{G}^2_b$ such that
\[
\mathbb{P} \bigl( g^{(0)}(X_T) \neq0 \bigr) \neq0.
\]

\item[$(\mathbf Ho)_2$] There exists some function $g
\in\mathcal{G}^3_b$ such that
\[
\mathbb{P} \bigl( g^\alpha(X_T) \neq0 \bigr) \neq0
\]
for $\alpha= (0)$, $(0,0)$ and $(j,0)$ for some $j \in\{1,\ldots,d\}$.
(Note that $g$ may be different for each $\alpha$.)

\item[$(\mathbf Ho)_3$] There exists some function $g
\in\mathcal{G}^4_b$ such that
\[
\mathbb{P} \bigl( g^\alpha(X_T) \neq0 \bigr) \neq0
\]
for $\alpha= (0)$, $(0,0)$ and $(j_1,0)$, $(j_2,0,0)$ for some
$(j_1,j_2) \in\{1,\ldots,d\}^2$. (Note that $g$ may be different for
each $\alpha$.)

Moreover, for any triplet $(\nu_1, \nu_2, \nu_3)\ne(0,0,0)$ we have
\[
\mathbb{P} \Biggl( \Biggl(\nu_1 g^{(0,0,0)} +\nu_2
f^y g^{(0,0) }+ \nu_3 \sum
_{\ell =1}^d f^{z^\ell} g^{(\ell,0,0) } \Biggr)
(X_T) \neq0 \Biggr) \neq0.
\]


\item[$(\mathbf Ho)_4$] There exists some function $g
\in\mathcal{G}^4_b$ such that
\[
\mathbb{P} \bigl( g^\alpha(X_T) \neq0 \bigr) \neq0
\]
for $\alpha= (0)$, $(0,0)$, $(j_1,0)$, $(j_2,0,0)$ for some $(j_1,j_2)
\in \{1,\ldots,d\}^2$. (Note that $g$~may be different for each
$\alpha$.)

Moreover, we have for pairs $(\nu_1,\nu_3)\ne(0,0)$, $(\nu_2, \nu
_4)\ne(0,0)$, 
\begin{eqnarray*}
&\displaystyle\mathbb{P} \Biggl( \Biggl(\nu_1 g^{(0,0,0)} +
\nu_3 \sum_{j=1}^d
{}^j v_g \Biggr) (X_T) \neq0 \Biggr)
\neq0,& 
\\
&\displaystyle\mathbb{P} \Biggl( \Biggl(\nu_2 g^{(\ell,0,0,0)} +
\nu_4 f^{z^\ell} \sum_{j=1}^d
g^{(j,0,0) } \Biggr) (X_T) \neq0 \Biggr) \neq0&
\end{eqnarray*}
for $1 \le\ell\le d$ and for any $(\nu_1$, $\nu_2$, $\nu_3$,
$\nu_4)\ne(0,0,0, 0)$ we have
\[
\mathbb{P} \Biggl( \Biggl(\nu_1 g^{(0,0,0,0)} + \nu_2
\sum_{j=1}^d {}^j
v_g^{(0)} + \nu_3 \sum
_{j=1}^d {}^j w_g +
\nu_4 \sum_{\ell=1}^d \sum
_{j=1}^d f^{z^\ell}{}^j
v_g^{(\ell)} \Biggr) (X_T) \neq0 \Biggr) \neq0,
\]
%
where we defined ${}^j v_g:= f^{z^j}g^{(j,0,0)}$ and ${}^j w_g:=
f^{z^j}g^{(j,0,0,0)}$, $1\le j \le d$.
\end{longlist}

%
%
\begin{Remark}
If the H\"ormander condition holds true, then all
conditions $(\mathbf Ho)_p$ are satisfied as the distribution of $X_T$ has
a smooth positive density with respect to the Lebesgue measure. 
\end{Remark}

\subsubsection{Description of the $H$-coefficients}\label{sec1.2.2}
We now specify the class of random variables $H$ used in the
Definition~\ref{de RK scheme} of the numerical schemes.

%
%
\begin{Definition}\label{de psi-H} \textup{(i)} For $m \ge0$, we denote by
$\mathcal{B}^{m}_{[0,1]}$ the set of bounded measurable functions
$\psi\dvtx[0,1] \rightarrow\mathbb{R}$ satisfying
\[
\int_{0}^{1}\psi(u)\,\mathrm{d}u = 1\quad
\mbox{and}\quad \mbox{if } m \ge1,\qquad \int_{0}^{1}
\psi(u)u^{k}\,\mathrm{d}u =0,\qquad 1\le k \le m.
\]

\textup{(ii)} Let $(\psi^{\ell})_{1\le\ell\le d} \in\mathcal
{B}^{m}_{[0,1]}$, for $t \in[0,T]$ and $h>0$ such that $t+h\le T$, we
define
\[
H^{\psi}_{t,h}:= \biggl(\frac1h\int_{t}^{t+h}
\psi^{\ell} \biggl(\frac
{u-t}{h} \biggr) \,\mathrm{d}W^\ell_{u}
\biggr)_{1\le\ell\le d},
\]
which is a row vector.

By convention, we set $H^{\psi}_{t,0} = 0$.
\end{Definition}

For a discussion on the choice of the above coefficients, we refer to
Remark~\ref{re justif H} and Section~\ref{sub se proxy Z}.
%

\subsubsection{One-stage schemes}\label{sec1.2.3}
We study here the order of the following family of schemes:
%
\begin{eqnarray*}
Y_{i } &=& \mathbb{E}_{t_{i }} \bigl[Y_{i +1} +
h_{i}b_1 f(Y_{i +1},Z_{i +1})
+h_{i} b_2f(Y_{i},Z_i) \bigr],
\\
Z_{i } &=& \mathbb{E}_{t_{i }} \bigl[H^{\psi_1}_{t_{i },h_{i}}
Y_{i +1} + h_{i}\beta_1 H^{\phi_1}_{t_{i },h_{i}}f(
Y_{i +1},Z_{i +1}) \bigr],
\end{eqnarray*}
where $\psi_1, \phi_1 \in\mathcal{B}^0_{[0,1]}$.

%
%
\begin{Theorem}\label{th s1}
\textup{(i)} Assume that $(\mathbf Hr)_1$ holds and that $\psi_1,
\phi_1 \in\mathcal{B}^0_{[0,1]}$. For $|\pi|$ small enough, the above
scheme is at least of order $1$ if
\[
1 = b_1 + b_2.
\]

Moreover, under $(\mathbf Ho)_1$, this condition is also necessary.

\textup{(ii)} Assume that $(\mathbf Hr)_1$ holds and that $\psi _1\in
\mathcal{B}^1_{[0,1]}$, $\phi_1 \in\mathcal{B}^0_{[0,1]}$. For $|\pi|$
small enough, the above scheme is at least of order $2$ if
\[
b_1 = b_2 = \tfrac12 \quad\mbox{and}\quad
\beta_1 = 1.
\]

Moreover, under $(\mathbf Ho)_2$, this condition is also necessary.
\end{Theorem}

%
%
\begin{Corollary} The above conditions lead to the following tableaux:
\[
\begin{array} {c@{\hspace*{5pt}}|@{\hspace*{5pt}}c@{
\hspace*{5pt}}|@{\hspace*{5pt}}c} 0 &0&0
\\
&&
\\[-15pt]
\hline1&1 & * \end{array} %
\quad\mbox{and}\quad%
\begin{array} {c@{\hspace*{5pt}}|@{\hspace*{5pt}}c@{\quad}c@{\hspace*{5pt}}|@{\hspace*{5pt}}c}
0 &0&0&0
\\
&&&
\\[-15pt]
\hline1&0&1 & * \end{array} %
\]
for the explicit Euler scheme and, respectively, the implicit version
and to the tableau
\[
\begin{array} {c@{\hspace*{5pt}}|@{\hspace*{5pt}}c@{\quad}c@{\hspace*{5pt}}|@{\hspace*{5pt}}c} 0 &0&0&0
\\
&&&
\\[-15pt]
\hline1&\frac12&\frac12 & 1 \end{array} %
\]
for the Crank--Nicholson scheme.
\end{Corollary}

%
\begin{Remark}
(i) The case of the Euler scheme has been widely studied in the
literature. Generally speaking, as soon as $f$ is Lipschitz-continuous,
the method has been shown to be convergent. Under weak regularity
assumption on the coefficient $g$, the order $\frac12$ can be
retrieved; see, for example,
\cite{boutou04,goblem05,gobmak10,zha04,geigob12}. The order~$1$
convergence has been first proved in \cite{goblab07} for the general
case when $f$ depends on~$Z$; see the references therein for the other
cases.

(ii) The Crank--Nicholson scheme of step (ii) has been studied in the
general case in \cite{criman10}. It is proved there to be of order $2$.

(iii) To the best of our knowledge, the necessary parts contained in
Theorem~\ref{th s1} are new.
\end{Remark}

%
\subsubsection{Two-stage schemes}\label{sec1.2.4}
We analyze here the order of the following family of schemes:

%
%
\begin{Definition}\label{de scheme explicit RK2}
%
\begin{eqnarray*}
Y_{i,2} &=& \mathbb{E}_{t_{i,2}} \bigl[Y_{i +1} +
a_{21} h_{i} f(Y_{i +1}, Z_{i +1})
\bigr]+ a_{22}h_{i}f(Y_{i,2},Z_{i,2}),
\\
Z_{i,2} &=& \mathbb{E}_{t_{i,2} } \bigl[H^{\psi_2}_{t_{i,2},c_2h_{i}}
Y_{i +1} + c_{2} h_i H^{\phi_2}_{t_{i,2},c_2h_{i}}f(Y_{i +1},
Z_{i
+1}) \bigr] 
\end{eqnarray*}
and
\begin{eqnarray*}
Y_{i } &=& \mathbb{E}_{t_{i }} \bigl[Y_{i } +
h_{i}b_1 f(Y_{i +1},Z_{i +1}) +
h_{i} b_2 f(Y_{i,2},Z_{i,2})
\bigr]+h_ib_3f(Y_{i},Z_i),
\\
Z_{i } &=& \mathbb{E}_{t_{i }} \bigl[H^{\psi_3}_{t_{i },h_{i}}
Y_{i +1} + \beta_1 H^{\phi_3}_{t_{i },h_{i}}h_if(
Y_{i +1},Z_{i +1}) + \beta_2H^{\phi_3}_{t_{i },(1-c_2)h_{i}}h_if(
Y_{i,2},Z_{i,2}) \bigr],
\end{eqnarray*}
where $\phi_2, \phi_3, \psi_2, \psi_3 \in\mathcal{B}^0$.
\end{Definition}


The following results concern implicit schemes (for the $Y$ part).

%
%
\begin{Theorem}
\label{th 2 stage scheme imp} \textup{(i)} Assume that $(\mathbf Hr)_3$
holds, $\psi_2,\psi_3 \in\mathcal{B}^2_{[0,1]}$, $\phi_2, \phi_3 \in
\mathcal{B} ^1_{[0,1]}$, $f^z = 0$ and $c_2 < 1$. For $|\pi|$ small
enough, the following conditions are sufficient to obtain at least an
order $3$ scheme
\begin{eqnarray*}
b_1&=&\frac12- \frac1{6 c_2},\qquad b_2 =
\frac1{6 c_2(1-c_2)},\qquad b_3 =\frac12-
\frac1{6 (1-c_2)},
\\
a_{21} &=& \frac{c_2}{2},\qquad \beta_{1} = 1-
\frac{1}{2 c_{2}},\qquad \beta_2 = \frac{1}{2 c_{2}}.
\end{eqnarray*}

\textup{(ii)} If, moreover, $(\mathbf Ho)_3$ holds, these conditions are also
necessary.\vadjust{\goodbreak}\enlargethispage{3pt}

\textup{(iii)} (Implicit order barrier) If $f^z\neq0$ and $(\mathbf Ho)_3$
holds, there is no order $3$ methods in the class of the schemes given
in Definition~\ref{de RK scheme} with
only two stages.\vspace*{-2pt}
\end{Theorem}

%
%
\begin{Corollary} 
\textup{(i)} For $0<c_2<1$, the above conditions lead to the following
tableau:
\[
%
\begin{array} {c@{\hspace*{5pt}}|@{\hspace*{5pt}}c@{\quad}c@{\quad}c@{\hspace*{5pt}}|@{\hspace*{5pt}}c@{\quad}c} 0 & & & & * &
\\
c_2& \frac{c_2}{2} & \frac{c_2}{2}\vspace*{2pt} &0 &
c_2 & *
\\[5pt]
&&&&&
\\[-20pt]
\hline
1 & \frac12-\frac1{6 c_2} & \frac1{6 c_2(1-c_2)}
& \frac{2-3c_2}{6(1-c_2)}& 1-\frac1{2c_2} & \frac1{2c_2}
\end{array}.
\]

\textup{(ii)} Observe that if $c_2=\frac23$, then $b_3 = 0$ and the tableau
has the following explicit form:
\[
\begin{array} {c@{\hspace*{5pt}}|@{\hspace*{5pt}}c@{\quad}c@{\hspace*{5pt}}|@{\hspace*{5pt}}c@{\quad}c} 0 & 0 & 0 & * & 0
\\
\frac23 & \frac13 & \frac13 & \frac23\vspace*{1pt} & *
\\[6pt]
&&&&
\\[-19pt]
\hline
\\[-12pt]
1 & \frac14& \frac34 & \frac14 & \frac34 \end{array}.  \vspace*{-2pt}
\]
\end{Corollary}

Part (iii) of the last theorem tells us that it is
generally not possible to get an order $3$ scheme with a two-stage
scheme, even if it is implicit, as soon as we have $f^z \neq0$. This
result differs from the ODE case. This fact is not surprising since the
schemes we consider are always explicit for the $Z$ part. The explicit
feature of the scheme and the related error, somehow propagates through
$f^z$. This will also be the case for schemes with a higher number of
stages. Since we are particularly interested in BSDEs with general
drivers, we see then that there is no advantage in using implicit
scheme instead of explicit ones. As a result, we concentrate from now
on in studying explicit schemes only.

The next result concerns then explicit schemes and
exhibits the similarity with the ODEs framework.\vspace*{-2pt}

%
%
\begin{Theorem}\label{th 2 stage scheme exp}
\textup{(i)} Assume that $(\mathbf Hr)_2$ holds and $\psi_2, \psi_3
\in\mathcal{B}^1_{[0,1]}$, $\phi_2$, \mbox{$\phi_3 \in
\mathcal{B}^0_{[0,1]}$}.

The scheme given in Definition~\ref{de RK scheme} is at least of order
2 if
\begin{eqnarray*}
&\displaystyle b_{1} = 1-\frac{1}{2 c_{2}} \quad\mbox{and}\quad
b_2 = \frac{1}{2c_{2}},&
\\[-2pt]
&\displaystyle \beta_1 + \beta_2\mathbf{1}_{\{c_2<1\}}
= 1.&
\end{eqnarray*}
%

\textup{(ii)} Moreover, if $(\mathbf Ho)_2$ holds, then the above
conditions are necessary.\vspace*{-2pt}
\end{Theorem}

It is easily checked that the above conditions leads to the following
tableau: For $0<c_2\le1$,
\[
%
\begin{array} {c@{\hspace*{5pt}}|@{\hspace*{5pt}}c@{\quad}c@{\hspace*{5pt}}|@{\hspace*{5pt}}c@{\quad}c} 0 & 0 & 0 & 0
& 0
\\
c_2& c_2 & 0 & c_2 & *
\\
&&&
\\[-14pt]
\hline
1 & 1-\frac1{2 c_2} & \frac1{2 c_2} &
\beta_1 & 1- \beta_1 \end{array} %
\]
with $\beta_1=1$ if $c_2=1$.\vadjust{\goodbreak}
%


\subsubsection{Three-stage schemes}\label{sec1.2.5}
We analyze next the order of the following family of schemes:

%
%
\begin{Definition}
%
%
\begin{eqnarray}
\label{eq RK3 explicit} Y_{i,2} &=& \mathbb{E}_{t_{i,2}}
\bigl[Y_{i +1} + h_{i}c_2 f(Y_{i +1},Z_{i +1})
\bigr],
\\
Z_{i,2} &=& \mathbb{E}_{t_{i,2}} \bigl[H^{\psi_2}_{t_{i,2},c_2h_{i}}
Y_{i +1} + h_{i}c_2 H^{\phi_2}_{t_{i,2},c_2h_{i}}
f( Y_{i +1},Z_{i +1}) \bigr],
\\
\label{eq RK3 explicit ti3} Y_{i,3} &=& \mathbb{E}_{t_{i,3}}
\bigl[Y_{i +1} + h_{i}a_{31} f(Y_{i +1},Z_{i,k})
+ h_{i}a_{32} f(Y_{i,2},Z_{i,2})
\bigr],
\\
Z_{i,3} &=& \mathbb{E}_{t_{i,3}} \bigl[H^{\psi_3}_{t_{i,3},c_3h_{i}}
Y_{i +1} + h_{i}\alpha_{31} H^{\phi_3}_{t_{i,3},c_3h_{i}}
f(Y_{i +1},Z_{i +1})
\nonumber
\\[-8pt]
\\[-8pt]
&&\hspace*{68pt} {} + h_{i}\tilde \alpha_{32}H^{\phi_3}_{t_{i,3},(c_3-c_2)h_{i}}
f(Y_{i,2},Z_{i,2}) \bigr].
\nonumber
\end{eqnarray}
The approximation at step~($\mathrm{i}$) is given by
%
%
\begin{eqnarray}
\label{eq RK3 explicit ti28} \qquad Y_{i } &=& \mathbb{E}_{t_{i }}
\bigl[Y_{i+1} + h_{i} \bigl( b_{1}f(Y_{i +1},Z_{i +1})
\nonumber
\\[-8pt]
\\[-8pt]
&&{}\hspace*{62pt} + b_{2}f( Y_{i,2},Z_{i,2})+
b_{3}f(Y_{i,3},Z_{i,3}) \bigr) \bigr],
\nonumber
\\
Z_{i } &=& \mathbb{E}_{t_{i }} \bigl[H^{\psi_4}_{t_{i },h_{i}}
Y_{i +1}
\nonumber
\\
&&\hspace*{15pt} {} + h_{i} \bigl( \beta_{1}
H^{\phi_4}_{t_{i },h_{i}} f(Y_{i +1},Z_{i +1})+ \tilde
\beta_{2} H^{\phi_4}_{t_{i },(1-c_2)h_{i}} f(Y_{i,2},Z_{i,2})
\\
&&\hspace*{140.5pt} {} + \tilde\beta_{3}H^{\phi_4}_{t_{i },(1-c_3)h_{i}}
f(Y_{i,3},Z_{i,3}) \bigr) \bigr]
\nonumber
\end{eqnarray}
with $\psi_2,\psi_3,\psi_4 \in B^2_{[0,1]}$, $\phi_2,\phi_3,\phi _4 \in B^1_{[0,1]}$.
\end{Definition}


%

%

%
%
\begin{Theorem}
\label{th 3 stage scheme}
\textup{(i)} Assume that $(\mathbf Hr)_3$ holds. The scheme given in
Definition~\ref{de RK scheme} is at least of order 3 if $c_2 \neq1$,
$c_2 \neq c_3$, and the following conditions hold true:
%
\begin{eqnarray*}
b_1+b_2+b_3 &=& 1,\qquad
b_2c_2 + b_3c_3 = \tfrac12,
\\
b_2c_2^2 + b_3c_3^2
&=& \tfrac13,\qquad b_3a_{32}c_2 =
b_3 \alpha_{32}c_2 = \tfrac16
\end{eqnarray*}
%
and
\begin{eqnarray*}
\beta_1 + \beta_2 + \beta_3
\mathbf{1}_{\{c_3 <1\}} &=& 1,
\\
\beta_2c_2 + \beta_3c_3
\mathbf{1}_{\{c_3 < 1\}} &=& \tfrac12.
\end{eqnarray*}

\textup{(ii)} Moreover, if $(\mathbf Ho)_3$ holds, 
then the above conditions are necessary.
\end{Theorem}

%
%
\begin{Remark}
\textup{(i)} If $c_2 = 1$, then $c_3=1$ and $\tilde\beta_2=
\tilde\beta_3 = 0$. Thus the approximation for $Z$ reads
\begin{eqnarray*}
Z_{i } &=& \mathbb{E}_{t_{i }} \bigl[H^{\psi_4}_{t_{i },h_{i}}
Y_{i +1} + h_{i}\beta_{1} H^{\phi_4}_{t_{i },h_{i}}
f( Y_{i +1},Z_{i +1}) \bigr].
\end{eqnarray*}
As shown in last section, this approximation leads generally to an
order 2 scheme only setting $\beta_1 = 1$.

\textup{(ii)} If $c_3 = c_2$, we obtain an order 2 scheme only as
well.\vadjust{\goodbreak}
\end{Remark}

Using \cite{but08} we get that
%
%
\begin{Corollary}\label{co
3 stage}
\textup{(i)} 
Assume that $c_2 \neq\frac23$, $c_3 \notin\{c_2,\frac23,1\}$. Then
the above conditions lead to the following tableau:
\[
\begin{array} {c@{\hspace*{5pt}}|@{\hspace*{5pt}}c@{\quad}c@{\quad}c} 0 & 0 & 0 & 0
\\
c_2& c_2 & 0 & 0
\\
c_3 & \frac{c_3(3c_2-3c_2^2-c_3)}{c_2(2-3c_2)} & \frac{c_3(c_3-c_2)}{c_2(2-3c_2)}\vspace*{3pt} & 0
\\[-5pt]
&&&
\\[-10pt]
\hline
1 & \frac{-3c_3+6 c_2 c_3 + 2 -3c_2}{6 c_2 c_3} & \frac{3 c_3- 2}{6 c_2 (c_3 - c_2)} & \frac{2 - 3c_2}{6c_3 (c_3-c_2)} \end{array},
\]
\[
\begin{array} {c@{\hspace*{5pt}}|@{\hspace*{5pt}}c@{\quad}c@{\quad}c} 0 & 0 & 0 & 0
\\
c_2& c_2 & * & 0
\\
c_3 &\frac{c_3(3c_2-3c_2^2-c_3)}{c_2(2-3c_2)}& \frac{c_3(c_3-c_2)}{c_2(2-3c_2)}\vspace*{3pt} &*
\\[-5pt]
&&&
\\[-10pt]
\hline
1 & \beta_1 & \frac{2c_3 - 1}{2 (c_3-c_2)}- \frac{c_3}{c_3-c_2}
\beta_1& \frac{c_3(1-2c_2)}{2c_3(c_3-c_2)} + \frac{c_2}{c_3-c_2} \beta_1
\end{array}.
\]


\textup{(ii)} If $c_3=1$ and $c_2 \neq\frac23$, then the above
conditions lead to the following tableau:
\[
\begin{array} {c@{\hspace*{5pt}}|@{\hspace*{5pt}}c@{\quad}c@{\quad}c@{\hspace*{5pt}}|@{\hspace*{5pt}}c@{\quad}c@{\quad}c} 0 & 0 & 0 & 0 & 0 & 0 & 0
\\
c_2& c_2 & 0 & 0 & c_2 & * & 0
\\
1 & \frac{(3c_2-3c_2^2-1)}{c_2(2-3c_2)} & \frac{1-c_2}{c_2(2-3c_2)} & 0& \frac{(3c_2-3c_2^2-1)}{c_2(2-3c_2)}&
\frac{1-c_2}{c_2(2-3c_2)}\vspace*{3pt} &*
\\[-5pt]
&&&&&
\\[-10pt]
\hline
1 & \frac{6 c_2 -3c_2 -1}{6 c_2 } & \frac{1}{6 c_2 (1 - c_2)} & \frac
{2 - 3c_2}{6 (1-c_2)} & 1-
\frac{1}{2c_2} &\frac{ 1}{2 c_2} & * \end{array}.
\]
\end{Corollary}


\subsubsection{Order barriers}\label{sec1.2.6}
As shown in the last sections, it is possible to derive explicit
methods of order $p=1,2,3$ using, respectively, $s=1,2,3$ stages. These
methods are optimal in the sense that $s<p$ is generally not possible
and $s>p$ would lead to more computational effort.

In the ODEs framework, such a result is well known; see \cite{but08}.
In fact, it is also known that it is possible to build explicit order 4
method using 4-stage schemes. A very interesting feature of explicit
methods is that to retrieve an order $p$ scheme with $p$ strictly
greater than $4$, one needs to use $s>p$ stages. This last result is
known as ``explicit order barriers''; see, for example, Theorem~370B in
\cite{but08}. Because ODEs are a special case of BSDEs, the same
explicit barriers will be encountered for BSDEs.

This leaves open the case $s=p=4$ for BSDEs. Theorem
\ref{th explicit barrier} below shows that generally $s > p$ already
for $p=4$ in the BSDEs framework. This means that the explicit barrier
is encountered earlier for BSDEs than for ODEs.

Before stating the main result of this section, let us also recall
part~(iii) of Theorem~\ref{th 2 stage scheme imp}, which reveals an
implicit order barrier in the BSDEs framework.

%
%
\begin{Proposition}[(Implicit barrier)] \label{pr implicit barrier}
Assume $(\mathbf Hr)_3$ holds and $f^z \neq0$, then there is no
implicit order 3 two-stage scheme, under the nondegeneracy
assumption~$(\mathbf Ho)_3$.\vadjust{\goodbreak}
\end{Proposition}



%
%
\begin{Theorem}[(Explicit barrier)] \label{th explicit barrier} We
assume that $f^y=0$ and $f^z \neq0$. There is no explicit four stage
methods in the class of methods given in Definition~\ref{de RK scheme}
which is of order $4$, provided that $(\mathbf Hr)_4$, $(\mathbf Ho)_4$
hold and that the $H$-coefficients are given by $H^i_j:=
H^{\psi_j}_{t_{i,j},c_j h_{i}}$ and $H^i_{j,k}:=
H^{\phi_j}_{t_{i,j},(c_j-c_k) h_{i}} $ with $\psi_j \in
\mathcal{B}^3_{[0,1]}$, $\phi_j \in\mathcal{B}^2_{[0,1]}$, $2 \le j
\le5$.
\end{Theorem}

%
%
\begin{Remark}\label{re explicit barrier}
Theorem~\ref{th explicit barrier} can be extended to the case of
$f^y\neq0$ and \mbox{$f^z \neq0$}. Indeed, the fact that $f^y \neq0$
will add more constraints to the problem. Note, however, that $(\mathbf
Ho)_4$ would need to be reformulated in this case.
\end{Remark}

\subsection{Outline}\label{sec1.3}

%
The rest of the paper is organized as follows. In Section~\ref{sec2},
we present some preliminary results used to study the order of
convergence. We also interpret the approximation of $Z$ as the
approximation of a proxy for $Z$ in dimension $d=1$.
Sections~\ref{sec3}--\ref{sec5} deal then with the proof of the order
for scheme with 1, 2~and~3 stages. Section~\ref{sec6} is dedicated to
the case of the four-stage methods and the proof of Theorem~\ref{th
explicit barrier}. Finally, the \hyperref[sec7]{Appendix} contains the proofs of the
results in Section~\ref{subse ge res} and the proofs of the preliminary
results.

\subsection{Notation}\label{sec1.4}

In the sequel $C$ is a positive constant whose value may change from
line to line depending on $T$, $d$, $\Lambda$, $X_0$ but which does not
depend on the choice of the partition $\pi$. We write $C_p$ if it
depends on some extra positive parameters $p$.

For $t \in\pi$, $R$ a random variable and $r$ a real
number, the notation $R=O_t(r)$ means that $|R| \le\lambda^\pi_t r$
where $\lambda^\pi_t$ is a positive random variable satisfying
\[
\mathbb{E} \bigl[ \bigl|\lambda^\pi_t \bigr|^p
\bigr] \le C_p
\]
for all $p > 0$, $t \in\pi$ and all partitions $\pi$.

The continuous and adapted process $U$ belongs to
$S^2([0,T])$ if
\[
\mathbb{E} \Bigl[\sup_{s \in[0,T]}|U_s|^2
\Bigr] < \infty.
\]

\subsubsection*{Multiple It\^o Integrals} For any process $U$ in
$S^2([0,T])$, we consider the following iterated Lebesgue--It\^o
integrals for a multi-index $\alpha$ with length $l$:
\[
I^\alpha_{t,s}[U]:= \cases{ U_s, &\quad if $l
=0$, \vspace*{3pt}
\cr
\displaystyle\int_t^s
I^{\alpha-}_{t,r}[U] \,\mathrm{d}r, &\quad if $l \ge1$ and $
\alpha_l = 0$, \vspace*{3pt}
\cr
\displaystyle\int
_t^s I^{\alpha-}_{t,r}[U] \,
\mathrm{d}W^j_r, &\quad if $l \ge1$ and $
\alpha_l = j$, $1 \le j \le d$.}
\]

One can recursively check that these integrals are well defined 
and that $I^\alpha[I^\beta[\cdot]] =I^{\beta*\alpha}[\cdot]$. We will
denote by $I^\alpha_{t,r}$ the multiple It\^o Integrals of the constant
process equal to one.

\subsection*{Abbreviation}
For $t \in[0,T]$, we denote $v^\alpha(t,X_t)$ by $v^\alpha_t$ and
$f^y(Y_t,Z_t)$ by $f^y_t$, where $f^y$ is the partial derivatives of
$f$ with respect to the variable $y$. Similarly $f^z_t:= f^z(Y_t,Z_t)$
where $f^z$ is the partial derivative of $f$ with respect to $z$.


\section{Preliminaries}\label{sec2}

%
%
%
%
%
%


\subsection{It\^o--Taylor expansions}\label{sec2.1}

The following proposition is Theorem 5.5.1 in \cite{klopla92} adapted
to our context.

%
\begin{Proposition}\label{pr basic expansion} Let $\mathcal{A}$ be a
hierarchical set and $\mathcal{B}(\mathcal{A})$ the associated
remainder set, for a function $v$ belonging to $\mathcal{G}^\beta_b$
for all $\beta\in \mathcal{B}(\mathcal{A})$. Then
\[
v(t+h,X_{t+h}) = \sum_{\alpha\in\mathcal{A}}v^\alpha_tI^{\alpha}_{t,t+h}
+ \sum_{\beta\in\mathcal{B}(\mathcal{A})}I^{\beta}_{t,t+h}
\bigl[v^\beta \bigr].
\]
\end{Proposition}

This leads to the following weak expansion formula:
%
%
\begin{Proposition}\label{pr weak expansion Y}
Let $m\ge0$. Then for a function $v \in\mathcal{G}^{m+1}_b$,
%
\[
\mathbb{E}_{t} \bigl[v(t+h,X_{t+h}) \bigr] 
= v_t +
hv_t^{(0)} + \frac{h^2}{2}v_t^{(0,0)}
+ \cdots+ \frac
{h^m}{m!}v_t^{(0)_m} + O_t
\bigl(h^{m+1} \bigr).
\]
\end{Proposition}

We now state another key expansion for the results below based on
Proposition~\ref{pr basic expansion} and Definition~\ref{de psi-H}.

%
%
\begin{Proposition}\label{pr weak expansion Z}
\textup{(i)} Let $m \ge0$, for $\psi=(\psi^\ell)_{1 \le\ell\le d}$ with
$\psi^\ell\in\mathcal{B}^m_{[0,1]}$, assuming that $v \in\mathcal
{G}^{m+2}_b$, then
%
\begin{eqnarray*}
\mathbb{E}_{t} \bigl[ \bigl(H^{\psi}_{t, h}
\bigr)^\ell v(t+ h,X_{t+ h}) \bigr] &=& v^{(\ell)}_t
+ h v^{(\ell,0)}_t + \cdots+ \frac{h^{m}}{m!}v^{(\ell)*(0)_{m}}_t
+ O_t \bigl(h^{m+1} \bigr).
\end{eqnarray*}

\textup{(ii)} For $\psi=(\psi^\ell)_{1 \le\ell\le d}$ with
$\psi^\ell\in \mathcal{B}^0_{[0,1]}$, assuming that $v
\in\mathcal{G}^1_b$, we have
\[
\mathbb{E}_{t} \bigl[ \bigl(H^{\psi}_{t, h}
\bigr)^\ell v(t+ h,X_{t+ h}) \bigr] = O_t(1).
\]

\textup{(iii)} If $L^{(0)}\circ L^{(\ell)} = L^{(\ell)}\circ L^{(0)}$,
for $\ell \in\{1,\ldots,d\}$, then the expansion of \textup{(i)} holds
true for $\psi= (1,\ldots,1) $.
\end{Proposition}

The proof of this proposition is postponed to the \hyperref[sec7]{Appendix}.

%
%
\begin{Remark}\label{re justif H}
\textup{(i)} The expansion of Proposition~\ref{pr weak expansion Z}(i)
motivates the definition of the $H$-coefficient. Indeed, we will apply
it to the functions $u$ and $u^{(0)}$ and are able to cancel the low
order term for a good choice of coefficients $(\alpha_{kj})$,
$(\beta_j)$; see the computations of the next sections.

\textup{(ii)} It is worth noticing that in the (very special) case
where $L^{(0)}\circ L^{(\ell)} = L^{(\ell)}\circ L^{(0)}$ for $\ell\in
\{1,\ldots,d\}$, one only needs to use in the definition of the scheme,
\mbox{$H$-}coefficients built with the function $\psi=(1,\ldots,1)$.
\end{Remark}

We conclude this paragraph by giving some examples of
function $\psi$ (\mbox{$d=1$}).

%
%
\begin{Example}\label{ex psi Bm}
\textup{(i)} The function $\psi=\mathbf{1}_{[0,1]}$ belongs to
$\mathcal {B}^{0}_{[0,1]}$.

\textup{(ii.a)} The polynomial function $x \mapsto\psi(x) = 4-6x$
belongs to $\mathcal{B}^{1}_{[0,1]}$.

\textup{(ii.b)} For $c \in(0,1)$, the function $\psi=
\frac1{c(c-1)}\mathbf{1}_{[1-c,1]} + \frac{c-2}{c-1}\mathbf
{1}_{[0,1]}$ belongs to~$\mathcal{B}^{1}_{[0,1]}$.
\end{Example}

\textup{(iii)} For $c, c' \in(0,1)$, $c \neq c'$,
%
\begin{eqnarray*}
\psi&=& \frac{1-c'}{c(1-c)(c'-c)}\mathbf{1}_{[1-c,1]} + \frac{c-1}{c' (1-c')(c'-c)}
\mathbf{1}_{[1-c',1]}
\\
&&{}+ \biggl(1+\frac{1}{(1-c)} + \frac {1}{(1-c')} \biggr)
\mathbf{1}_{[0,1]}
\end{eqnarray*}
belongs to $\mathcal{B}^{2}_{[0,1]}$.

\subsection{A class of proxy for $Z$}\label{sub se proxy Z}

The solution of the BSDE (\ref{eq bsde YZ}) consists in the pair
process $(Y,Z)$. Unlike $Y$, the second component is not ``directly
available'' in~(\ref{eq bsde YZ}) since it is defined as the integrand
in the martingale part. However, we can use~(\ref{eq bsde YZ}) to
construct first a proxy for $Z$. As we shall see, the sequence of
processes $(Z_i)_{i\le n}$ are discrete-time approximation of this
proxy. The results below are based on the expansion given in
Proposition~\ref{pr weak expansion Z}. The discussion in this section
assumes $d=1$.

%
%
\begin{Definition}\label{de proxy} For $m \ge0$, let $\psi\in
\mathcal{B}^{m}_{[0,1]}$
%
%
\begin{equation}
\label{eq de proxy} Z_{t,h}^{\psi}:=\mathbb{E}_{t}
\biggl[H^{\psi}_{t,h}\int_{t}^{t+h}
Z_{u} \,\mathrm{d}W_{u} \biggr].\vadjust{\goodbreak}
\end{equation}
\end{Definition}

For later use, we denote $H^{\psi}_{t,h}(u) =
\mathbb{E}_{u} [H^{\psi}_{t,h} ]$, $t\le u \le t+h$.

%
%
\begin{Proposition}\label{pr proxy}
Let $m \ge0$, and assume that $u \in\mathcal{G}^{m+2}_b$. For $\psi
\in
\mathcal{B}^{m}_{[0,1]}$, the following holds:
\[
Z_{t} = Z_{t,h}^{\psi} + O \bigl(h^{m+1}
\bigr).
\]
\end{Proposition}

\begin{pf} One observes that
\[
Z_{t,h}^{\psi} =\frac1h \mathbb{E}_{t} \biggl[\int
_{t}^{t+h} \psi \biggl(\frac{s-t}{h} \biggr)
Z_{s}\,\mathrm{d}s \biggr] = \frac{1}{h} \mathbb{E}_{t}
\biggl[\int_{t}^{t+h} \psi \biggl(
\frac{s-t}{h} \biggr) u^{(1)}_{s}\,\mathrm{d}s \biggr].\vadjust{\goodbreak}
\]

Applying the expansion given in Proposition~\ref{pr weak expansion Y}
to $u^{(1)}$ up to order $m$ and using the assumption on $\psi$, we
obtain
\begin{eqnarray*}
Z_{t,h}^{\psi} &=& \sum_{k=0}^{m}
u^{(0)_k * (1)}_{t} \frac1h \int_{t}^{t+h}
\psi \biggl(\frac{s-t}{h} \biggr)\frac{(s-t)^{k}}{k!} \,\mathrm {d}s +
O_{t} \bigl(h^{m+1} \bigr)
\\
&=& \sum_{k=0}^{m} u^{(0)_k * (1)}_{t}
\frac{h^{k}}{k!} \int_{0}^{1}
\psi(r)r^{k} \,\mathrm{d}r + O_{t} \bigl(h^{m+1}
\bigr)
\\
&=&Z_{t} + O_{t} \bigl(|h|^{m+1} \bigr),
\end{eqnarray*}
recalling that $Z_t = u^{(1)}_t$.
\end{pf}

%
%
\begin{Remark}\label{re todo proxy Z} 
Of course one can build other types of proxies for $Z$ based
on~(\ref{eq de proxy}), for example, at $t=0$,
\[
\mathbb{E} \biggl[H^\psi_{0,h}\int_0^h
Z_s \,\mathrm{d}W_s+ \lambda_1
Z_h + \lambda_2 Z_{h/2} \cdots \biggr].
\]
In this case, $\psi$ will be required to satisfy different constraints
in order to obtain the desired order of convergence.
\end{Remark}
%

It remains to derive the discrete-time approximation $(Z_i)$.

Observe that, using (\ref{eq bsde YZ}),
%
%
\begin{eqnarray}
Z^\psi_{t,h} &:=& \mathbb{E}_{t}
\biggl[H^{\psi}_{t,h}\int_{t}^{t+h}
Z_{u} \,\mathrm{d} W_{u} \biggr]
\nonumber
\\[-8pt]
\label{eq show begins}
\\[-8pt]
&=& \mathbb{E}_{t} \biggl[H^{\psi
}_{t,h} \biggl(
Y_{t+h} + \int_{t}^{t+h}f(Y_{u},Z_{u})
\,\mathrm{d}u \biggr) \biggr].
\nonumber
\end{eqnarray}

In \cite{boutou04,goblab07}, the approximation of the $Z$ process is
given by
\[
\bar{Z}^{\mathbf{1}}_{t_{i },h_{i}}:= \mathbb{E}_{t_{i }}
\bigl[H^{\mathbf{1}}_{t_{i },h_{i}}Y_{t_{i +1}} \bigr].
\]

In order to obtain high-order approximation of the process $Z$, we
discretize the integral term in the right-hand side in (\ref{eq show
begins}), with $t = t_{i }$. For $\psi\in\mathcal {B}^m_{[0,1]}$, $m
\ge1$, we will approximate this term by the following:
%
%
\begin{equation}
\label{eq de approx Hf } h\sum_{j=1}^q
\beta_j \mathbb{E}_{t_{i }} \bigl[H^{\phi_j}_{t_{i }, (1-c_j) h_{i} }
f(Y_{t_{i,j}},Z_{t_{i,j}}) \bigr],
\end{equation}
where the coefficients $\beta_j \in\mathbb{R}$ and the function
$\phi_j$ belongs to $\mathcal{B}^{m-1}_{[0,1]}$, for $1 \le j \le q$.

%
%
\begin{Remark}
Alternatively, one can approximate directly
\[
\mathbb{E}_{t} \biggl[H^{\psi}_{t_{i },h_{i}}\int
_{t_{i }}^{t_{i +1}}f(Y_u,Z_u) \,
\mathrm{d} u \biggr] = \mathbb{E}_{t_{i }} \biggl[\int
_{t_{i }}^{t_{i +1}}H^{\psi}_{t_{i
},h_{i} }(u)f(Y_u,Z_u)
\,\mathrm{d}u \biggr]
\]
by
\[
h\sum_{j=1}^q \beta_j
\mathbb{E}_{t_{i }} \bigl[H^{\psi}_{t_{i },h_{i}}(t_{i,j})
f(Y_{t_{i,j}},Z_{t_{i,j}}) \bigr].
\]
However, since generally $H^{\psi}_{t_{i },h_{i}}(t_{i,j}) \neq
H^{\psi}_{t_{i },(1-c_j)h_{i}}$, one would then require stronger
assumptions on the function $\psi$ and the $H$-coefficient which, in
turn, will lead to higher computational complexity.
\end{Remark}

The approximation given in (\ref{eq de approx Hf }) is still
theoretical since it uses the true value $Y_{t_{i,j}}$ and
$Z_{t_{i,j}}$. We need to introduce several stages to obtain
approximations of these intermediate values.



\section{One-stage schemes}\label{sec3}

\subsection{\texorpdfstring{Proof of Theorem \protect\ref{th s1}\textup{(i)}}
{Proof of Theorem 1.3(i)}}\label{sec3.1}

\textup{(1)} We first compute the error expansion for the $Z$ part of
the scheme. By (\ref{eq de hY hZ}), we have, for $1 \le\ell\le d$,
\begin{eqnarray*}
\hat{Z}_{t_{i }}^\ell&:=& \mathbb{E}_{t_{i }} \bigl[
\bigl(H^{\psi_1}_{t_{i },h_{i}} \bigr)^\ell{Y}_{t_{i +1}} +
h_{i}\beta_1 \bigl(H^{\phi_1}_{t_{i },h_{i}}
\bigr)^\ell f({Y}_{t_{i +1}},{Z}_{t_{i +1}}) \bigr]
\\
&=& \mathbb{E}_{t_{i }} \bigl[ \bigl(H^{\psi_1}_{t_{i },h_{i}}
\bigr)^\ell u_{t_{i +1}} - h_{i} \beta_1
\bigl(H^{\phi_1}_{t_{i },h_{i}} \bigr)^\ell u^{(0)}_{t_{i +1}}
\bigr],
\end{eqnarray*}
recalling (\ref{b_pde_t}).

Using Proposition~\ref{pr weak expansion Z}, we get
%
%
\begin{equation}
\label{eq o1-1s Z1} \hat{Z}_{t_{i }}^\ell =
Z_{t_{i }}^\ell+ O_{t_{i
}} \bigl(|\pi| \bigr),
\end{equation}
since $u \in\mathcal{G}^2_b$, recalling $(\mathbf Hr)_1$ and $\psi
^\ell_1, \phi^\ell_1 \in\mathcal{B}^{0}_{[0,1]}$.

This basically means that as soon as $\psi^\ell_1
\in\mathcal{B}^0_{[0,1]}$, $1 \le\ell\le d$, the choice of $\beta
_1$~is arbitrary. Indeed, by definition of the truncation\vadjust{\goodbreak} error for the
$Z$ component [see (\ref{eq de trunc error loc YnZ})--(\ref{eq de trunc
error YnZ})], we have
\[
\mathcal{T}_Z(\pi) = O \bigl(|\pi| \bigr),
\]
which is the order we aim to obtain.

\textup{(2a)} We now compute the error expansion for the $Y$-part.
First observe that
\begin{eqnarray*}
\hat{Y}_{t_{i }} &:=& \mathbb{E}_{t_{i }} \bigl[{Y}_{t_{i +1}}
+ h_{i}b_1 f({Y}_{t_{i +1}},{Z}_{t_{i +1}})
+h_{i}b_2f(\hat{Y}_{t_{i }},\hat{Z}_{t_{i }})
\bigr]
\\
&=& \mathbb{E}_{t_{i }} \bigl[{Y}_{t_{i +1}} + h_{i}b_1
f({Y}_{t_{i +1}},{Z}_{t_{i +1}}) +h_{i} b_2f({Y}_{t_{i }},
\hat{Z}_{t_{i }}) \bigr] +h_{i}b_2\delta
f_{t_{i }}, 
\end{eqnarray*}
where $\delta f_{t_{i }} =f(\hat{Y}_{t_{i }},\hat{Z}_{t_{i }})
-f(Y_{t_{i }},\hat{Z}_{t_{i }})$.
This leads to
\[
\hat{Y}_{t_{i }}:= \mathbb{E}_{t_{i }} \bigl[u_{t_{i +1}} -
h_{i}b_1 u^{(0)}_{t_{i +1}}-
h_{i}b_2u^{(0)}_{t_{i }} \bigr]
+h_{i}b_2\delta f_{t_{i }}.
\]

Using Proposition~\ref{pr weak expansion Y}, we compute
\begin{eqnarray*}
\hat{Y}_{t_{i }} &:=& u_{t_{i }} +h_{i}(1-
b_1- b_2) u^{(0)}_{t_{i }} 
+h_{i}b_2\delta f_{t_{i }} +
O_{t_{i }} \bigl(|\pi|^2 \bigr).
\end{eqnarray*}
Since $f$ is Lipschitz-continuous and $u^{(0)}$ bounded, we obtain for
$|\pi|$ small enough that $\hat{Y}_{t_{i }} = {Y}_{t_{i }} + O_{t_{i
}}(|\pi|)$ which implies that $\delta f_{t_{i }} = O_{t_{i }}(|\pi|)$
and thus
%
%
\begin{equation}
\label{eq s1o1 Y} \hat{Y}_{t_{i }}:= Y_{t_{i }}
+h_{i}(1- b_1- b_2) u^{(0)}_{t_{i }}
+ O_{t_{i }} \bigl(|\pi|^2 \bigr).
\end{equation}
The condition $b_1+b_2 =1$ is thus sufficient to retrieve at least an
order-$1$ scheme.

\textup{(2b)} Under $(\mathbf Ho)_1$, this condition is also necessary.

Indeed, combining definition (\ref{eq de trunc error loc YnZ})--(\ref
{eq de trunc error YnZ}) and (\ref{eq s1o1 Y}), we compute
\[
\mathcal{T}_Y(\pi) = \sum_{i=0}^{n-1}
h_{i}|1- b_1- b_2|^2 \mathbb{E}
\bigl[ \bigl|u^{(0)}_{t_{i }} \bigr|^2 \bigr] + O \bigl(|
\pi|^2 \bigr).
\]
Interpreting the sum in the last equation as a Riemann sum and taking
the limit as $|\pi|\rightarrow0$,
we obtain
\[
\lim_{|\pi|\downarrow0} \mathcal{T}_Y(\pi) = |1-
b_1- b_2|^2 \int_0^T
\mathbb{E} \bigl[ \bigl|u^{(0)}(t,X_t) \bigr|^2 \bigr]
\, \mathrm{d}t.
\]

If $ (1- b_1- b_2)^2\neq0$, since the scheme must be of order 1, we
must have
\[
\int_0^T \mathbb{E} \bigl[
\bigl|u^{(0)}(t,X_t) \bigr|^2 \bigr]\,\mathrm{d}t = 0
\]
for solutions $u$ of (\ref{b_pde_t}) such that $u \in\mathcal{G}^2_b$,
recalling Definition~\ref{de order}. In particular, at $t=T$, since $t
\mapsto\mathbb{E} [|u^{(0)}(t,X_t)|^2 ]$ is continuous, we get
\[
\mathbb{E} \bigl[ \bigl|g^{(0)}(X_T) \bigr|^2 \bigr] = 0
\qquad\mbox{for all } g \in\mathcal{G}^2_b.
\]
Under $(\mathbf Ho)_1$, this yields a contradiction.

\subsection{\texorpdfstring{Proof of Theorem \protect\ref{th s1}\textup{(ii)}}
{Proof of Theorem 1.3(ii)}}\label{sec3.2}

\textup{(1a)} We first compute the expansion for the $Z$~part. By
definition [see (\ref{eq de hY hZ})], we have, for $1 \le\ell\le d$,
\begin{eqnarray*}
\hat{Z}_{t_{i }}^\ell&:=& \mathbb{E}_{t_{i }} \bigl[
\bigl(H^{\psi_1}_{t_{i },h_{i}} \bigr)^\ell{Y}_{t_{i +1}} +
h_{i}\beta_1 \bigl(H^{\phi_1}_{t_{i },h_{i}}
\bigr)^\ell f({Y}_{t_{i +1}},{Z}_{t_{i +1}}) \bigr]
\\
&=& \mathbb{E}_{t_{i }} \bigl[ \bigl(H^{\psi_1}_{t_{i },h_{i}}
\bigr)^\ell u_{t_{i +1}} - h_{i} \beta_1
\bigl(H^{\phi_1}_{t_{i },h_{i}} \bigr)^\ell u^{(0)}_{t_{i +1}}
\bigr].
\end{eqnarray*}
Using Proposition~\ref{pr weak expansion Z}, we have
%
%
\begin{eqnarray}
\label{eq o21s Z1} \hat{Z}_{t_{i }}^\ell&=&
Z_{t_{i }}^\ell+h_{i}(1- \beta_1)
u^{(\ell,0)}_{t_{i }} + O_{t_{i }} \bigl(| \pi|^2
\bigr)
\end{eqnarray}
since $u \in\mathcal{G}^3_b$ and $\psi^\ell_1 \in\mathcal{B}^1_{[0,1]}$,
$\phi^\ell_1 \in\mathcal{B}^0_{[0,1]}$, $1 \le\ell\le d$.\vadjust{\goodbreak}

Using a first-order Taylor expansion, this leads to
%
%
\begin{eqnarray}
\label{eq o21s Z2} f({Y}_{t_{i }},\hat{Z}_{t_{i }}) &=& -
u^{(0)}_{t_{i }} + h_{i}(1-\beta_1)\sum
_{\ell
= 1}^d f^{z^\ell}_{t_{i }}
u^{(\ell,0)}_{t_{i }} + O_{t_{i }} \bigl(|\pi|^2
\bigr),
\end{eqnarray}
recalling that $f \in C^2_b$ under $(\mathbf Hr)_2$.

From (\ref{eq o21s Z1}) we deduce that the condition $1-\beta_1=0$ is
sufficient to obtain $\mathcal{T}_Z(\pi)= O(|\pi|^2)$, recalling
(\ref{eq de trunc error loc YnZ})--(\ref{eq de trunc error YnZ}).

\textup{(1b)} If we assume that $(\mathbf Ho)_2$ holds, this condition
is also necessary. Indeed, one computes that
\begin{eqnarray*}
\frac{\mathcal{T}_Z(\pi)}{|\pi|^2} &=& \sum_{i=0}^{n-1}h_{i}(1-
\beta_1)^2 \sum_{\ell=1}^d
\mathbb{E} \bigl[ \bigl|u^{(\ell,0)}_{t_{i }} \bigr|^2 \bigr] + O
\bigl(|\pi|^2 \bigr)
\end{eqnarray*}
for grids with constant mesh size.

Then by interpreting the sum in the last equation as a Riemann sum, we
obtain
\begin{eqnarray*}
\lim_{|\pi| \downarrow0} \frac{\mathcal{T}_Z(\pi)}{|\pi|^2} &=& (1-\beta_1)^2
\int_0^T \sum_{\ell=1}^d
\mathbb{E} \bigl[ \bigl|u^{(\ell,0)}(t,X_t) \bigr|^2 \bigr]
\, \mathrm{d}t,
\end{eqnarray*}
where the limit is taken over the grids with constant mesh size. If $
(1- \beta_1)^2\neq0$, since we are looking at a scheme of order 2, we
must have
\begin{eqnarray*}
\int_0^T \sum_{\ell=1}^d
\mathbb{E} \bigl[ \bigl|u^{(\ell,0)}(t,X_t) \bigr|^2 \bigr]
\, \mathrm{d}t &=& 0
\end{eqnarray*}
for the solution $u$ of (\ref{b_pde_t}) such that $u
\in\mathcal{G}^3_b$, recalling Definition~\ref{de order}. In
particular, at $t=T$, since $t \mapsto\sum_{\ell=1}^d \mathbb{E}
[|u^{(\ell,0)}(t,X_t)|^2 ]$ is continuous, we get
\[
\sum_{\ell=1}^d \mathbb{E} \bigl[
\bigl|g^{(\ell,0)}(X_T) \bigr|^2 \bigr] = 0\qquad\mbox{for all
} g \in\mathcal{G}^3_b.
\]
Under $(\mathbf Ho)_2$, this yields a contradiction.


We assume now that the condition $\beta_1=1$ holds.

\textup{(2a)} For the $Y$-part, we have
\begin{eqnarray*}
\hat{Y}_{t_{i }} &:=& \mathbb{E}_{t_{i }} \bigl[{Y}_{t_{i +1}}
+ h_{i}b_1 f({Y}_{t_{i +1}},{Z}_{t_{i +1}})
+h_{i}b_2f(\hat{Y}_{t_{i }},\hat{Z}_{t_{i }})
\bigr]
\\
&=& \mathbb{E}_{t_{i }} \bigl[{Y}_{t_{i +1}} + h_{i}b_1
f({Y}_{t_{i +1}},{Z}_{t_{i +1}}) +h_{i} b_2f({Y}_{t_{i }},
\hat{Z}_{t_{i }}) \bigr] +h_{i}b_2\delta
f_{t_{i }}, 
\end{eqnarray*}
where $\delta f_{t_{i }} =f(\hat{Y}_{t_{i }},\hat{Z}_{t_{i }})
-f(Y_{t_{i }},\hat{Z}_{t_{i }})$.

Combining the last equality with (\ref{eq o21s Z2}) and recalling that
$\beta_1=1$, we get
\begin{eqnarray*}
\hat{Y}_{t_{i }} &:=& \mathbb{E}_{t_{i }} \bigl[u_{t_{i +1}} -
h_{i}b_1 u^{(0)}_{t_{i +1}}-
h_{i}b_2u^{(0)}_{t_{i }} \bigr]
+h_{i}b_2\delta f_{t_{i }} +
O_{t_{i }} \bigl(|\pi|^3 \bigr).
\end{eqnarray*}
Since $u \in\mathcal{G}^3_b$, we use Proposition~\ref{pr weak expansion
Y} to compute
%
%
\begin{eqnarray}
\hat{Y}_{t_{i }} &=& {Y}_{t_{i }} + h_{i}(1-b_1-b_2)u^{(0)}_{t_{i }}
+ h_{i}^2 \bigl(\tfrac12-b_1
\bigr)u^{(0,0)}_{t_{i }} 
\nonumber
\\[-8pt]
\label{eq o21s Y1}
\\[-8pt]
&&{} +h_{i}b_2 \delta f_{t_{i }} +
O_{t_{i }} \bigl(|\pi|^3 \bigr).
\nonumber
\end{eqnarray}
%
We observe that $\hat{Y}_{t_{i }} = {Y}_{t_{i }} + O(|\pi|)$ which
leads to
\[
\delta f_{t_{i }} = O_{t_{i }} \bigl(|\pi| \bigr)
\]
since $f$ is Lipschitz-continuous.

Combining (\ref{eq o21s Y1}) with the last estimate, we obtain
\[
\hat{Y}_{t_{i }} = {Y}_{t_{i }} + h_{i}(1-b_1-b_2)u^{(0)}_{t_{i }}
+ O_{t_{i }} \bigl(|\pi|^2 \bigr).
\]
The condition
\[
(1-b_1-b_2) = 0
\]
is sufficient to obtain a method at least of order $1$.

\textup{(2b)} Using the same arguments as in step (2b) of the proof of
part (i) of Theorem~\ref{th s1}, we obtain that this condition is
necessary if
$(\mathbf Ho)_2$ holds. 

\textup{(2c)} We thus assume from now on that this condition holds, and
we get
%
%
\begin{equation}
\label{eq o21s Y2} \hat{Y}_{t_{i }} = {Y}_{t_{i }} +
O_{t_{i }} \bigl(| \pi|^2 \bigr),
\end{equation}
which leads, since $f$ is Lipschitz-continuous, to 
$ \delta f_{t_{i }} = O_{t_{i }}(|\pi|^2). $ Inserting this estimate
back into
(\ref{eq o21s Y1}), we obtain
%
%
\begin{eqnarray}
\label{eq o21s Y3} \hat{Y}_{t_{i }} &=& {Y}_{t_{i }} +
h_{i}^2 \bigl(\tfrac12-b_1 \bigr)u^{(0,0)}_{t_{i }}
+ O_{t_{i }} \bigl(|\pi|^3 \bigr),
\end{eqnarray}
recalling that $b_1+b_2 = 1$.

The condition $\frac12-b_1= 0$ is therefore sufficient to obtain a
method at least of order $2$.

\textup{(2d)} If we assume that $(\mathbf Ho)_2$ holds, this condition
is also necessary. Indeed, one computes that
\begin{eqnarray*}
\frac{\mathcal{T}_Y(\pi)}{|\pi|} &=& \sum_{i=0}^{n-1}h_{i}
\biggl(\frac12-b_1 \biggr)^2 \mathbb{E} \bigl[
\bigl|u^{(0,0)}_{t_{i }} \bigr|^2 \bigr] + O \bigl(|
\pi|^2 \bigr)
\end{eqnarray*}
for grids $\pi$ with constant mesh size.\vadjust{\goodbreak}

Then, as the limit of a Riemann sum, we obtain that
\begin{eqnarray*}
\lim_{|\pi| \downarrow0} \frac{\mathcal{T}_Y(\pi)}{|\pi|} &=& \biggl(\frac
12-b_1 \biggr)^2\int_0^T
\mathbb{E} \bigl[ \bigl|u^{(0,0)}(t,X_t) \bigr|^2 \bigr]
\, \mathrm{d}t,
\end{eqnarray*}
where the limit is taken over the grids with constant mesh size. If
$\frac12-b_1\neq0$, since the scheme must be of order 2, we must have
\[
\int_0^T \mathbb{E} \bigl[
\bigl|u^{(0,0)}(t,X_t) \bigr|^2 \bigr]\,\mathrm{d}t = 0
\]
for solution u of (\ref{b_pde_t}) such that $u \in\mathcal{G}^3_b$,
recalling Definition~\ref{de order}. In particular, at $t=T$, since $t
\mapsto\mathbb{E} [|u^{(0,0)}(t,X_t)|^2 ]$ is continuous, we get
\[
\mathbb{E} \bigl[ \bigl|g^{(0,0)}(X_T) \bigr|^2 \bigr] = 0
\qquad\mbox{for all } g \in\mathcal{G}^3_b.
\]
Under $(\mathbf Ho)_2$, this yields a contradiction and completes the proof
of the theorem. 

\section{Two-stage schemes}\label{sec4}


\subsection{\texorpdfstring{Proof of Theorem \protect\ref{th 2 stage scheme imp}}
{Proof of Theorem 1.4}}\label{sec4.1}

(1a) We first compute the error expansion at the intermediary step
(step $j=2$), recalling that $(\mathbf Hr)_3$ is in force.

For $1 \le\ell\le d$, we have that
\begin{eqnarray*}
\hat{Z}_{t_{i,2}}^\ell&:=& \mathbb{E}_{t_{i,2}} \bigl[
\bigl(H^{\psi
_2}_{t_{i,2},c_2h_{i} } \bigr)^\ell Y_{t_{i +1}} +
h_{i}c_2 \bigl(H^{\phi
_2}_{t_{i,2},c_2h_{i}}
\bigr)^\ell f(Y_{t_{i +1}},Z_{t_{i +1}}) \bigr]
\\
& =& \mathbb{E}_{t_{i,2}} \bigl[ \bigl(H^{\psi_2}_{t_{i,2},c_2h_{i}}
\bigr)^\ell u_{t_{i +1}} - h_{i}c_2
\bigl(H^{\phi_2}_{t_{i,2},c_2h_{i}} \bigr)^\ell u^{(0)}_{t_{i +1}}
\bigr].
\end{eqnarray*}
Since $u\in\mathcal{G}^4_b$, we apply Proposition~\ref{pr weak
expansion Z} and get, for $1 \le\ell\le d$,
%
%
\begin{eqnarray}
\label{eq 2s exp Zti2} \hat{Z}_{t_{i,2}}^\ell
&=& \frac{c_2^2}{2}h_{i}^2
u^{(\ell,0,0)}_{t_{i,2}} + O_{t_{i,2}} \bigl(|\pi|^3
\bigr).
\end{eqnarray}
Using a first order Taylor expansion, we obtain 
\begin{eqnarray*}
f(Y_{t_{i,2}}, \hat{Z}_{t_{i,2}} ) & =& -u^{(0)}_{t_{i,2}}
- \frac{c_2^2}{2}h_{i}^2 \sum
_{\ell= 1}^d f^{z^\ell}_{t_{i,2}}u^{(\ell,0,0)}_{t_{i,2}}+
O_{t_{i,2}} \bigl(|\pi|^3 \bigr),
\end{eqnarray*}
recalling that $f \in C^2_b$.

(1b) For the $Y$-part, we have, denoting $\delta f_{t_{i,2}}
=f(\hat{Y}_{t_{i,2}},\hat{Z}_{t_{i,2}})
-f(Y_{t_{i,2}},\hat{Z}_{t_{i,2}}) $,
\begin{eqnarray*}
\hat{Y}_{t_{i,2}} &:=& \mathbb{E}_{t_{i,2}} \bigl[Y_{t_{i +1}} +
h_{i}a_{21} f(Y_{t_{i +1}},Z_{t_{i +1}})
\bigr]+a_{22} h_{i} f(Y_{t_{i,2}},
\hat{Z}_{t_{i,2}})+a_{22} h_{i} \delta
f_{t_{i,2}}
\\
&=& \mathbb{E}_{t_{i,2}} \bigl[u_{t_{i +1}} - h_{i} \bigl(
a_{21}u^{(0)}_{t_{i +1}} + a_{22}
u^{(0)}_{t_{i,2}} \bigr) \bigr] 
+a_{22} h_{i} \delta f_{t_{i,2}} + O_{t_{i,2}}
\bigl(|\pi|^3 \bigr).
\end{eqnarray*}

Using Proposition~\ref{pr weak expansion Y}, we compute
%
%
\begin{eqnarray}
\label{eq interm s2 Y} \hat{Y}_{t_{i,2}} &=&Y_{t_{i,2}} + \biggl(
\frac{c_2^2}{2}-a_{21}c_2 \biggr)h_{i}
^2u^{(0,0)}_{t_{i,2}} 
+a_{22} h_{i} \delta f_{t_{i,2}} +
O_{t_{i,2}} \bigl(|\pi|^3 \bigr),
\end{eqnarray}
recalling that $u \in\mathcal{G}^3_b$.

Since $f$ is Lipschitz continuous, we get that $\delta f_{t_{i,2}} =
O_{t_{i,2}}(|\pi|^2)$.

Inserting this estimate back into (\ref{eq interm s2 Y}), we obtain
\begin{eqnarray*}
\hat{Y}_{t_{i,2}} &=&Y_{t_{i 2}} + \biggl(\frac{c_2^2}{2}-a_{21}c_2
\biggr)h_{i} ^2u^{(0,0)}_{t_{i,2}}
+ O_{t_{i,2}} \bigl(|
\pi|^3 \bigr).
\end{eqnarray*}

Combining a first-order Taylor expansion with the last equality and
(\ref{eq 2s exp Zti2}) leads~to
%
%
\begin{eqnarray}
f(\hat{Y}_{t_{i,2}},\hat{Z}_{t_{i,2}}) &=& -u^{(0)}_{t_{i,2}}
+ \biggl(\frac
{c_2^2}{2}-a_{21}c_2
\biggr)h_{i}^2 f^y_{t_{i 2}}u^{(0,0)}_{t_{i,2}}
\nonumber
\\[-8pt]
\label{eq fti2}
\\[-8pt]
&&{}- \frac{c_2^2}{2}h_{i}^2 \sum
_{\ell=1}^df^{z^\ell}_{t_{i 2}}u^{(\ell,0,0)}_{t_{i,2}}
+ O_{t_{i,2}} \bigl(|\pi|^3 \bigr).
\nonumber
\end{eqnarray}

(2a) We now study the error at the final step for the $Z$-part.


We compute the following expansion, for $1 \le\ell\le d$:
%
\begin{eqnarray*}
\hat{Z}_{t_{i }}^\ell&:=& \mathbb{E}_{t_{i }} \bigl[
\bigl(H^{\psi_3}_{t_{i },h_{i}} \bigr)^\ell Y_{t_{i +1}}+
\beta_1 h_{i} \bigl(H^{\phi_3}_{t_{i },h_{i}}
\bigr)^\ell f(Y_{t_{i +1}},Z_{t_{i +1}})
\\
&&\hspace*{57pt} {} + \tilde{ \beta}_2 h_{i}
\bigl(H^{\phi_3}_{t_{i },(1-c_2)h_{i}} \bigr)^\ell f(
\hat{Y}_{t_{i,2}},\hat{Z}_{t_{i,2}} ) \bigr]
\\
&=& \mathbb{E}_{t_{i }} \bigl[ \bigl(H^{\psi_3}_{t_{i },h_{i}}
\bigr)^\ell u_{t_{i +1}} - \beta_1 h_{i}
\bigl(H^{\phi_3}_{t_{i },h_{i}} \bigr)^\ell u^{(0)}_{t_{i +1}}
- \tilde{\beta}_2 h_{i} \bigl(H^{\phi_3}_{t_{i },(1-c_2)h_{i}}
\bigr)^\ell u^{(0)}_{t_{i,2}} \bigr]
\\
&&{}+ O_{t_{i }} \bigl(|\pi|^3 \bigr),
\end{eqnarray*}
where we used (\ref{eq fti2}), Proposition~\ref{pr weak expansion Z}
and $(\mathbf Hr)_3$, observing that $f^yu^{(0,0)}$ and
$f^{z^\ell}u^{(\ell,0,0)}$, $1 \le\ell\le d$, belong to
$\mathcal{G}^1_b$.

Using Proposition~\ref{pr weak expansion Z} again, we obtain, for $1
\le\ell\le d$,
%
%
\begin{eqnarray}
\hat{Z}_{t_{i }}^\ell-
Z_{t_{i }}^\ell&:=& (1-\beta_1 - \tilde
\beta_2)h_{i} u^{(\ell,0)}_{t_{i }}
\nonumber
\\[-8pt]\label{eq 2so2 exp Z}
\\[-8pt]
&&{} + \bigl(\tfrac12-\beta_1 - (1-c_2)\tilde
\beta_2 \bigr) h_{i}^2 u^{(\ell,0,0)}_{t_{i }}
+O_{t_{i }} \bigl(|\pi|^3 \bigr).
\nonumber
\end{eqnarray}
%

(2b) For the local truncation error on the $Z$-part to be of order $2$,
recalling (\ref{eq de trunc error YnZ}), it is clear, according to
(\ref{eq 2so2 exp Z}), that the following condition is sufficient:
%
%
\begin{equation}
\label{eq cond 2s Z o2} 1-\beta_1 - \tilde\beta_2 = 0.
\end{equation}
Similarly, to retrieve local truncation error on the $Z$-part to be of
order $3$, the following conditions are sufficient:
%
%
\begin{eqnarray}
1 -\beta_1 - \tilde\beta_2 &=& 0, \label{eq cond 2s Z
o3_1}
\\
c_2 \tilde\beta_2- \tfrac12&=& 0. \label{eq cond 2s Z
o3_2}
\end{eqnarray}

(2c) We now prove that condition (\ref{eq cond 2s Z o2}) is necessary
to obtain an order $2$ scheme under $(\mathbf Ho)_2$, recalling that
$(\mathbf Ho)_3$ implies $(\mathbf Ho)_2$. We compute, for grids with constant
mesh size,
\[
\frac{\mathcal{T}_Z(\pi)}{|\pi|^2} = (1-\beta_1 - \tilde\beta_2
)^2h_{i} \sum_{\ell= 1}^d
\mathbb{E} \bigl[ \bigl|u^{(\ell,0)}_{t_{i }} \bigr|^2 \bigr] +
O_{t_{i }} \bigl(|\pi|^2 \bigr)
\]
and then (Riemann sum)
\[
\lim_{|\pi| \downarrow0} \frac{\mathcal{T}_Z(\pi)}{|\pi|^2} = (1-\beta_1 -
\tilde\beta_2 )^2 \int_0^T
\sum_{\ell= 1}^d \mathbb{E} \bigl[
\bigl|u^{(\ell,0)}(t,X_t) \bigr|^2 \bigr]\,\mathrm{d}t.
\]
If $ (1-\beta_1 - \tilde\beta_2)^2\neq0$, since the scheme must be of
order 2, we must have
\[
\int_0^T \sum_{\ell= 1}^d
\mathbb{E} \bigl[ \bigl|u^{(\ell,0)}(t,X_t) \bigr|^2 \bigr]
\, \mathrm{d}t = 0.
\]
In particular, at $t=T$, since $t \mapsto\sum_{\ell= 1}^d
\mathbb{E} [|u^{(\ell,0)}(t,X_t)|^2 ]$ is continuous, we get
\[
\sum_{\ell= 1}^d \mathbb{E} \bigl[
\bigl|g^{(\ell,0)}(X_T) \bigr|^2 \bigr] = 0 \qquad \mbox{for
all } g \in\mathcal{G}^3_b.
\]
Under $(\mathbf Ho)_2$, this yields a contradiction. %

(2d) Under $(\mathbf Ho)_3$, it is thus necessary that
$(1-\beta_1 - \tilde\beta_2 )^2 = 0$ to retrieve an order $2$ and a
fortiori an order $3$ schemes. The expansion error for the $Z$ part
reads then
%
%
\begin{eqnarray}
\hat{Z}_{t_{i }}^\ell-
Z_{t_{i }}^\ell&:=& \bigl(c_2 \tilde
\beta_2- \tfrac12 \bigr) h_{i}^2
u^{(\ell,0,0)}_{t_{i }} +O_{t_{i }} \bigl(|\pi|^3
\bigr),\qquad 1 \le\ell\le d.
\end{eqnarray}

Using the same techniques as in step (2c), one will get that condition
(\ref{eq cond 2s Z o3_2}) is necessary to obtain an order $3$ scheme
under $(\mathbf Ho)_3$.
%

(3) We study the error expansion on the $Y$ part at the final step. We
aim to obtain an order $3$ scheme. From the definition of the
truncation error, it is obviously necessary that the local truncation
error for the $Z$ part is of order $3$. We work then under this
condition [see step (2d)] and then we have
%
%
\begin{equation}
\label{eq 2so2 exp Z bis} f(Y_{t_{i }},\hat{Z}_{t_{i }}) = -
u^{(0)}_{t_i} + O_{t_{i }} \bigl(|\pi|^3
\bigr).
\end{equation}

For the $Y$-part, using (\ref{eq fti2}) and (\ref{eq 2so2 exp Z bis}),
we have
\begin{eqnarray*}
\hat{Y}_{t_{i }} &:=& \mathbb{E}_{t_{i }} \bigl[Y_{t_{i +1}} +
h_{i}b_1 f(Y_{t_{i +1}},Z_{t_{i +1}})+
h_{i}b_2f(\hat{Y}_{t_{i,2}},\hat{Z}_{t_{i,2}}
) \bigr]
\\
&&{} +h_ib_{3}f({Y}_{t_{i }},\hat
{Z}_{t_{i }} ) + \delta f_{t_{i }}
\\
&=& \mathbb{E}_{t_{i }} \Biggl[u_{t_{i +1}} - h_{i}b_1
u^{(0)}_{t_{i +1}} - h_{i}b_1
u^{(0)}_{t_{i,2}}
\\
&&\hspace*{17pt} {} + \biggl(\frac{c_2^2}{2}-a_{21}c_2
\biggr)h_{i}^3 f^yu^{(0,0)}_{t_{i,2}}
- \frac
{c_2^2}{2}h_{i}^3 \sum
_{\ell=1}^d f^{z^\ell}u^{(\ell,0,0)}_{t_{i,2}}
\Biggr]
\\
&&{}-h_ib_{3} u^{(0)}_{t_i} +
h_{i}\delta f_{t_{i }} + O_{t_{i }} \bigl(|
\pi|^4 \bigr).
\end{eqnarray*}

Using Proposition~\ref{pr weak expansion Y} and $(\mathbf Hr)_3$,
\begin{eqnarray*}
\hat{Y}_{t_{i }} &=& \mathbb{E}_{t_{i }} \bigl[u_{t_{i +1}} -
h_{i}b_1 u^{(0)}_{t_{i +1}} -
h_{i}b_1 u^{(0)}_{t_{i,2}} -
h_ib_{3} u^{(0)}_{t_i} \bigr]
\\
&&{}+ \biggl(\frac{c_2^2}{2}-a_{21}c_2
\biggr)h_{i}^3 f^yu^{(0,0)}_{t_{i }}
- \frac
{c_2^2}{2}h_{i}^3 \sum
_{\ell=1}^d f^{z^\ell}u^{(\ell,0,0)}_{t_{i }}
+ h_{i}\delta f_{t_{i }} + O_{t_{i }} \bigl(|
\pi|^4 \bigr).
\end{eqnarray*}

Using Proposition~\ref{pr weak expansion Y}, we get
%
%
\begin{eqnarray}
\label{eq 2so2 Y exp} \qquad\hat{Y}_{t_{i }} - Y_{t_{i }} &=&
(1-b_1-b_2-b_3)h_{i}u^{(0)}_{t_{i }}
+ \biggl( \frac12-b_1-b_2(1-c_2) \biggr)
h_{i}^2 u^{(0,0)}_{t_{i }}
\nonumber
\\
&&{}+ \biggl( \frac16-\frac{b_1+b_2(1-c_2)^2}{2} \biggr) h_{i}^3
u^{(0,0,0)}_{t_{i }}
\nonumber
\\[-8pt]
\\[-8pt]
&&{} + \biggl(\frac{c_2^2}{2}-a_{21}c_2
\biggr)h_{i}^3 f^yu^{(0,0)}_{t_{i }}
\nonumber
\\
&&{}- \frac{c_2^2}{2}h_{i}^3 \sum
_{\ell=1}^d f^{z^\ell}u^{(\ell,0,0)}_{t_{i }}
+ h_{i}\delta f_{t_{i }} + O_{t_{i }} \bigl(|
\pi|^4 \bigr).
\nonumber
\end{eqnarray}

Using the last equation, we obtain that $\delta f_{t_{i }} =
O_{t_{i }}(|\pi|)$, which leads to
\[
\hat{Y}_{t_{i }} - Y_{t_{i }} = (1-b_1-b_2-b_3)h_{i}u^{(0)}_{t_{i }}
+ O_{t_{i }} \bigl(|\pi|^2 \bigr).
\]
Under $(\mathbf Ho)_3$, it appears then that the following condition is
necessary to retrieve an order $\ge1$ scheme:
%
%
\begin{equation}
\label{eq cond 2so3 1} b_1+b_2+b_3 = 1.
\end{equation}
We then assume that this condition holds and obtain
\[
\hat{Y}_{t_{i }} - Y_{t_{i }} =h_{i}\delta
f_{t_{i }} + O_{t_{i }} \bigl(|\pi|^2 \bigr).
\]
We thus compute
\[
\delta f_{t_{i }} = h_{i}f^y \delta
f_{t_{i }} + O_{t_{i }} \bigl(|\pi|^2 \bigr).
\]
And for $|\pi|$ small enough, $\delta f_{t_{i }} = O_{t_{i }}(|\pi|^2)$.
Inserting this into (\ref{eq 2so2 Y exp}) and recalling that (\ref{eq
cond 2so3 1}) is in force, we get that
\[
\hat{Y}_{t_{i }} - Y_{t_{i }} = \bigl(\tfrac12 - b_1 -
b_2(1-c_2) \bigr) h_{i}^2
u^{(0,0)}_{t_{i }} + O_{t_{i }} \bigl(|\pi|^3
\bigr).
\]

Under $(\mathbf Ho)_3$, the condition $\frac12 =b_1 + b_2(1-c_2) $ is then
necessary to obtain an order $2$ scheme,\vadjust{\goodbreak} and we thus assume it holds.
Arguing as before we now obtain $\delta f_{t_{i }} = O_{t_{i }}(|\pi
|^3)$ and
then
\begingroup
\abovedisplayskip=7pt
\belowdisplayskip=7pt
%
%
\begin{eqnarray}
\label{eq 2so2 Y exp bis} \qquad\quad\hat{Y}_{t_{i }} - Y_{t_{i }}
&=& \biggl(\frac{b_2(1-c_2)c_2}{2} - \frac1{12} \biggr) h_{i}^3
u^{(0,0,0)}_{t_{i }} + \biggl(\frac{c_2^2}{2}-a_{21}c_2
\biggr)h_{i}^3 f^yu^{(0,0)}_{t_{i }}
\nonumber
\\[-9pt]
\\[-9pt]
&&{}- \frac{c_2^2}{2}h_{i}^3 \sum
_{\ell=1}^d f^{z^\ell}u^{(\ell,0,0)}_{t_{i }}
+ O_{t_{i }} \bigl(|\pi|^4 \bigr).
\nonumber
\end{eqnarray}

(3b) If $f^{z^\ell} = 0$ for all $\ell\in\{1,\ldots, d\}$, one obtains
that
\[
b_2 = \frac1{6(1-c_2)c_2} \quad\mbox{and}
\quad a_{21} = \frac{c_2}{2}
\]
are sufficient conditions for the methods to be of order 3.

Under $(\mathbf Ho)_3$, these are also necessary conditions.

This completes the proof of (i) and (ii).

(4) To prove (iii), we use (\ref{eq 2so2 Y exp bis}) again. We observe
that under $(\mathbf Ho)_3$, if $f^{z^\ell} \neq0$ for some $\ell\in
\{1,\ldots,d\}$, since $c_2 >0$, the methods is at most of
order~$2$.\vspace*{-2pt}

\subsection{\texorpdfstring{Proof of Theorem \protect\ref{th 2 stage scheme exp}}
{Proof of Theorem 1.5}}\label{sec4.2}\vspace*{-2pt}

\mbox{}
\begin{pf} The computation for the explicit case is almost the
same---easier, in fact. The main difference comes from the fact that we
are only interested in order 2 schemes. We thus need a bit less
regularity. Following the step of the last proof, one then gets the
following error expansion:
%
%
\begin{eqnarray}
\hat{Y}_{t_{i }} - Y_{t_{i }} &=& (1-b_1-b_2)h_{i}u^{(0)}_{t_{i }}
\nonumber
\\[-9pt]
\\[-9pt]
&&{}+ \bigl(\tfrac12-b_1-b_2(1-c_2) \bigr)
h_{i}^2 u^{(0,0)}_{t_{i }} +
O_{t_{i }} \bigl(|\pi|^3 \bigr)
\nonumber
\end{eqnarray}
and, for $\ell\in\{1,\ldots,d\}$,
%
%
\begin{eqnarray}
\hat{Z}_{t_{i }}^\ell- Z_{t_{i }}^\ell&:=&
(1-\beta_1 - \tilde\beta_2)h_{i}
u^{(\ell,0)}_{t_{i }} 
+O_{t_{i }} \bigl(|\pi|^2 \bigr).
\end{eqnarray}
%
Under $(\mathbf Hr)_2$, the conditions
\[
1-b_1-b_2 = 0,\qquad b_2c_2 =
\tfrac12,\qquad 1-\beta_1 - \tilde \beta_2 =0
\]
%
are obviously sufficient. Under $(\mathbf Ho)_2$, using the same
techniques as in steps \mbox{(2c)--(2d)} of the proof of
Theorem~\ref{th 2 stage scheme imp}, one proves that these conditions
are necessary, which completes the proof of the theorem.\vspace*{-2pt}
\end{pf}

\section{Three-stage schemes}\label{sec5}\vspace*{-2pt}


\subsection{\texorpdfstring{Proof of Theorem \protect\ref{th 3 stage scheme}}
{Proof of Theorem 1.6}}\label{sec5.1}\vspace*{-2pt}


(1a) We compute the error expansion at the intermediary
step $j=2$.
\begin{eqnarray*}
\hat{Y}_{t_{i,2}} &:=& \mathbb{E}_{t_{i,2}} \bigl[Y_{t_{i +1}} +
h_{i}c_2 f(Y_{t_{i +1}},Z_{t_{i +1}}) \bigr] =
\mathbb{E}_{t_{i,2}} \bigl[u_{t_{i +1}} - h_{i}c_2
u^{(0)}_{t_{i +1}} \bigr],\vadjust{\goodbreak}
\\
{Z}_{t_{i,2}} &:=& \mathbb{E}_{t_{i,2}} \bigl[H^{\psi_2}_{t_{i,2},c_2h_{i}}
Y_{t_{i +1}}+ h_{i}c_2 H^{\psi_2}_{t_{i,2},c_2h_{i}}f(Y_{t_{i
+1}},Z_{t_{i +1}})
\bigr]
\\
& =& \mathbb{E}_{t_{i,2}} \bigl[H^{\psi_2}_{t_{i,2},c_2h_{i}}
u_{t_{i +1}} - h_{i} c_2 H^{\psi_2}_{t_{i,2},c_2h_{i}}
u^{(0)}_{t_{i +1}} \bigr].
\end{eqnarray*}\endgroup
Under $(\mathbf Hr)_3$, applying Propositions~\ref{pr weak expansion
Y}~and~\ref{pr weak expansion Z}, we have
\begin{eqnarray*}
\hat{Y}_{t_{i,2}} &=& {Y}_{t_{i,2}} - \frac
{c_2^2}{2}h_{i}^2u^{(0,0)}_{t_{i,2}}
+ O_{t_{i,2}} \bigl(|\pi|^3 \bigr),
\\
\hat{Z}_{t_{i,2}}^\ell&=& {Z}_{t_{i,2}}^\ell-
\frac
{c_2^2}{2}h_{i}^2u^{(\ell,0,0)}_{t_{i,2}} +
O_{t_{i,2}} \bigl(|\pi|^3 \bigr),\qquad\ell\in\{1,\ldots,d\},
\\
f(\hat{Y}_{t_{i,2}},\hat{Z}_{t_{i,2}}) &:=& -u^{(0)}_{t_{i,2}}
- \frac{c_2^2}{2}h_{i}^2 \Biggl(f^y
u^{(0,0)}_{t_{i,2}}+ \sum_{\ell= 1}^d
f^{z^\ell}_{t_{i,2}} u^{(\ell,0,0)}_{t_{i,2}} \Biggr)+
O_{t_{i,2}} \bigl(|\pi|^3 \bigr).
\end{eqnarray*}

(1b) Error expansion at step $3$.
\begin{eqnarray*}
\hat{Y}_{t_{i,3}} &:=& \mathbb{E}_{t_{i,3}} \bigl[Y_{t_{i +1}} +
h_{i}a_{31} f(Y_{t_{i +1}},Z_{t_{i +1}}) +
h_{i}a_{32}f(Y_{t_{i,2}},Z_{t_{i,2}}) \bigr],
\\
%
\hat{Z}_{t_{i,3}} &:=& \mathbb{E}_{t_{i,3}} \bigl[H^{\psi
_3}_{t_{i,3},c_3h_{i}}
Y_{t_{i +1}} + h_{i}\alpha_{31} H^{\phi_3}_{t_{i,3},c_3h_{i}}
f(Y_{t_{i +1}},Z_{t_{i +1}})
\\
&&\hspace*{66pt} {} + h_{i} \tilde \alpha_{32}
H^{\phi_3}_{t_{i,3},(c_3-c_2)h_{i}} f(Y_{t_{i,2}},Z_{t_{i,2}}) \bigr].
\end{eqnarray*}
With this definition and using step (1a), we compute
\begin{eqnarray*}
\hat{Y}_{t_{i,3}} &=& \mathbb{E}_{t_{i,3}} \bigl[u_{t_{i +1}} -
h_{i}a_{31} u^{(0)}_{t_{i +1}} -
h_{i}a_{32} u^{(0)}_{t_{i,2}} \bigr] +
O_{t_{i,3}} \bigl(|\pi|^3 \bigr),
\end{eqnarray*}
which leads to, recalling $a_{31}+a_{32} = c_3$,
\begin{eqnarray*}
\hat{Y}_{t_{i,3}} &=& Y_{t_{i,3}} - \biggl(\frac{c_3^2}{2}-c_2
a_{32} \biggr)h_{i} ^2u^{(0,0)}_{t_{i,3}}
+ O_{t_{i,3}} \bigl(|\pi|^3 \bigr).
\end{eqnarray*}

Equivalently, we get, for $\ell\in\{1,\ldots,d\}$,
\begin{eqnarray*}
\hat{Z}_{t_{i,3}}^\ell&=& Z_{t_{i,3}}^\ell-
\biggl(\frac{c_3^2}{2}-c_2 \tilde\alpha_{32}
\biggr)h_{i}^2u^{(\ell,0,0)}_{t_{i,3}} +
O_{t_{i,3}} \bigl(|\pi|^3 \bigr).
\end{eqnarray*}
And we obtain
\begin{eqnarray*}
f(\hat{Y}_{t_{i,3}},\hat{Z}_{t_{i,3}} )&:=& - u^{(0)}_{t_{i,3}}
- \biggl(\frac{c_3^2}{2}-c_2 a_{32}
\biggr)h_{i}^2 f^y_{t_{i,3}}
u^{(0,0)}_{t_{i,3}}
\\
&&{} - \biggl(\frac{c_3^2}{2}-c_2 \tilde\alpha_{32}
\biggr) h_{i}^2 \sum_{\ell= 1}^d
f^{z^\ell}_{t_{i,3}}u^{(\ell,0,0)}_{t_{i,3}}
\\
&&{} + O_{t_{i,3}} \bigl(|\pi|^3 \bigr).
\end{eqnarray*}

(1c) Error expansion at the final step for $Z$.
\begin{eqnarray*}
\hat{Y}_{t_{i }} &:=& \mathbb{E}_{t_{i }} \bigl[Y_{t_{i +1}} +
h_{i}b_1 f({Y}_{t_{i +1}},{Z}_{t_{i +1}} ) +
h_{i}b_{2} f(\hat{Y}_{t_{i,2}},\hat{Z}_{t_{i,2}}
)+ h_{i}b_{3} f(\hat{Y}_{t_{i,3}},\hat
{Z}_{t_{i,3}} ) \bigr],
\\
\hat{Z}_{t_{i }} &:=& \mathbb{E}_{t_{i }} \bigl[H^{\psi_4}_{t_{i },h_{i}}
Y_{t_{i +1}} + h_{i}\beta_{1}H^{\phi_4}_{t_{i },h_{i}}
f({Y}_{t_{i +1}},{Z}_{t_{i +1}} )
\\
&&\hspace*{15pt} {} + h_{i} \beta_{2} H^{\phi_4}_{t_{i },(1-c_2)h_{i}}
f(\hat{Y}_{t_{i,2}}, \hat{Z}_{t_{i,2}} ) + h_{i}\tilde
\beta_{3} H^{\phi_4}_{t_{i },(1-c_3)h_{i}} f(\hat{Y}_{t_{i,3}},
\hat{Z}_{t_{i,3}} ) \bigr].
\end{eqnarray*}

Using the results of step (1), we then compute, for
$\ell\in\{1,\ldots,\ell\}$,
\begin{eqnarray*}
\hat{Z}_{t_{i }}^\ell&=& \mathbb{E}_{t_{i }} \bigl[
\bigl(H^{\psi_4}_{t_{i },h_{i}} \bigr)^\ell u_{t_{i +1}}-
h_{i}\beta_{1} \bigl(H^{\phi_4}_{t_{i },h_{i}}
\bigr)^\ell u^{(0)}_{t_{i +1}}
\\
&&\hspace*{15pt} {} -h_{i} \beta_{2} \bigl(H^{\phi_4}_{t_{i
},(1-c_2)h_{i}}
\bigr)^\ell u^{(0)}_{t_{i,2}} - h_{i}\tilde
\beta_{3} \bigl(H^{\phi_4}_{t_{i },(1-c_3)h_{i}} \bigr)^\ell
u^{(0)}_{t_{i,3}} \bigr]
\\
&&{}- \beta_2 \frac{c_2^2}{2}h_{i}^3
\mathbb{E}_{t_{i }} \Biggl[ \bigl(H^{\phi
_4}_{t_{i },(1-c_2)h_{i}}
\bigr)^\ell \Biggl(f^y_{t_{i,2}}u^{(0,0)}_{t_{i,2}}+
\sum_{j=1}^d f^{z^j}_{t_{i,2}}
u^{(j,0,0)}_{t_{i,2}} \Biggr) \Biggr]
\\
&&{}- \tilde\beta_3 \biggl(\frac{c_3^2}{2}-c_2
a_{32} \biggr)h_{i}^3 \mathbb
{E}_{t_{i }} \bigl[ \bigl(H^{\phi_4}_{t_{i },(1-c_3)h_{i}}
\bigr)^\ell f^y_{t_{i,3}}u^{(0,0)}_{t_{i,3}}
\bigr]
\\
&&{} - \tilde\beta_3 \biggl(\frac{c_3^2}{2}-c_2
\tilde\alpha_{32} \biggr) h_{i}^3
\mathbb{E}_{t_{i }} \Biggl[ \bigl(H^{\phi_4}_{t_{i },(1-c_3)h_{i}}
\bigr)^\ell\sum_{j=1}^d
f^{z^j}_{t_{i,3}}u^{(j,0,0)}_{t_{i,3}} \Biggr]
\\
&&{}+O_{t_{i }} \bigl(|\pi|^3 \bigr).
\end{eqnarray*}
Under $(\mathbf Hr)_3$, since $f^yu^{(0,0)}, f^{z^j}u^{(j,0,0)} \in
\mathcal{G}^1_b$, $j \in\{1,\ldots,d\}$, we obtain using Proposition
\ref{pr weak expansion Z}(ii), for all $\ell\in\{1,\ldots,d\}$,
\begin{eqnarray*}
\mathbb{E}_{t_{i }} \Biggl[ \bigl(H^{\phi_4}_{t_{i },(1-c_2)h_{i}}
\bigr)^\ell \Biggl(f^y_{t_{i,2}}u^{(0,0)}_{t_{i,2}}+
\sum_{j=1}^d f^{z^j}_{t_{i,2}}
u^{(j,0,0)}_{t_{i,2}} \Biggr) \Biggr] &=& O_{t_{i }}(1),
\\
\mathbb{E}_{t_{i }} \bigl[ \bigl(H^{\phi_4}_{t_{i },(1-c_3)h_{i}}
\bigr)^\ell f^y_{t_{i,3}}u^{(0,0)}_{t_{i,3}}
\bigr] &=& O_{t_{i }}(1)
\end{eqnarray*}
and
\[
\mathbb{E}_{t_{i }} \Biggl[ \bigl(H^{\phi_4}_{t_{i
},(1-c_3)h_{i}}
\bigr)^\ell\sum_{j=1}^d
f^{z^j}_{t_{i,3}}u^{(j,0,0)}_{t_{i,3}} \Biggr] =
O_{t_{i }}(1).
\]


And then
\begin{eqnarray*}
\bar{Z}_{t_{i }} &=& \mathbb{E}_{t_{i }} \bigl[H^{\psi_4}_{t_{i },h_{i}}
u_{t_{i +1}} - h_{i}\beta_{1} H^{\phi_4}_{t_{i },h_{i}}
u^{(0)}_{t_{i +1}}
\\
&&\hspace*{15pt} {}- h_{i}\beta_{2}H^{\phi_4}_{t_{i },(1-c_2)h_{i}}
u^{(0)}_{t_{i,2}} - h_{i}\tilde\beta_{3}
H^{\phi_4}_{t_{i },(1-c_3)h_{i}} u^{(0)}_{t_{i,3}} \bigr]
\\
&&{} + O_{t_{i }} \bigl(| \pi|^3 \bigr).
\end{eqnarray*}

Using the expansion of Proposition~\ref{pr weak expansion Z}, this
leads to the following truncation error for the $Z$ part:
%
%
\begin{eqnarray}
\label{eq RK3 trunc error Z} \qquad\quad\mathcal{T}_Z(\pi)&:= &\sum
_i (1- \beta_1 + \beta_2 +
\tilde\beta_3)^2h_{i} ^3 \sum
_{\ell=1}^d \mathbb{E} \bigl[ \bigl|
u^{(\ell,0)}_{t_{i }} \bigr|^2 \bigr]
\nonumber
\\
&&{} +\sum_i \biggl(\frac12-\beta_1 -
\beta_2(1-c_2) - \tilde\beta_3(1-c_3)
\biggr)^2h_{i}^5\sum
_{\ell=1}^d \mathbb{E} \bigl[ \bigl| u^{(\ell,0,0)}_{t_{i }}
\bigr|^2 \bigr]
\\
&&{} + O \bigl(|\pi|^6 \bigr).
\nonumber
\end{eqnarray}


(1d) Error expansion at the final step for $Y$.
\[
\hat{Y}_{t_{i }}:= \mathbb{E}_{t_{i }} \bigl[Y_{t_{i +1}} +
h_{i}b_1 f({Y}_{t_{i +1}},{Z}_{t_{i +1}} ) +
h_{i}b_{2} f(\hat{Y}_{t_{i,2}},\hat{Z}_{t_{i,2}}
)+ h_{i}b_{3} f(\hat{Y}_{t_{i,3}},\hat
{Z}_{t_{i,3}} ) \bigr].
\]

We compute that
\begin{eqnarray*}
\hat{Y}_{t_{i }} &=& \mathbb{E}_{t_{i }} \bigl[u_{t_{i +1}} -
h_{i}b_{1} u^{(0)}_{t_{i +1}} -
h_{i}b_2 u^{(0)}_{t_{i,2}} -
h_{i}b_3 u^{(0)}_{t_{i,3}} \bigr]
\\
&&{}- b_2 \frac{c_2^2}{2}h_{i}^3
\mathbb{E}_{t_{i }} \Biggl[ \Biggl(f^y_{t_{i,2}}u^{(0,0)}_{t_{i,2}}+
\sum_{\ell=1}^df^{z^\ell}_{t_{i,2}}u^{(\ell,0,0)}_{t_{i,2}}
\Biggr) \Biggr]
\\
&&{}- b_3 \biggl(\frac{c_3^2}{2}-c_2
a_{32} \biggr)h_{i}^3 \mathbb{E}_{t_{i }}
\bigl[f^y_{t_{i,3}} u^{(0,0)}_{t_{i,3}} \bigr]
\\
&&{} - b_3 \biggl(\frac{c_3^2}{2}-c_2 \tilde
\alpha_{32} \biggr)h_{i}^3\mathbb{E}_{t_{i }}
\Biggl[\sum_{\ell=1}^df^{z^\ell}_{t_{i,3}}
u^{(\ell,0,0)}_{t_{i,3}} \Biggr]
\\
&&{}+O_{t_{i }} \bigl(|\pi|^4 \bigr),
\end{eqnarray*}
which leads to
\begin{eqnarray*}
\hat{Y}_{t_{i }}&=& \mathbb{E}_{t_{i }} \bigl[u_{t_{i +1}} -
h_{i}b_{1} u^{(0)}_{t_{i +1}} -
h_{i}b_2 u^{(0)}_{t_{i,2}} -
h_{i}b_3 u^{(0)}_{t_{i,3}} \bigr]
\\
&&{}- \biggl(b_2 \frac{c_2^2}{2} + b_3
\frac{c_3^2}{2} - b_3 c_2 a_{32}
\biggr)h_{i}^3 f^y u^{(0,0)}_{t_{i }}
\\
&&{}- \biggl(b_2 \frac{c_2^2}{2} + b_3
\frac{c_3^2}{2} - b_3 c_2 \tilde\alpha
_{32} \biggr)h_{i}^3 \sum
_{\ell=1}^df^{z^\ell}_{t_{i }}
u^{(\ell,0,0)}_{t_{i }}
\\
&&{}+O_{t_{i }} \bigl(|\pi|^4 \bigr).
\end{eqnarray*}

Using then Proposition~\ref{pr weak expansion Y}, we obtain the
following global truncation error for~$Y$:
%
%
\begin{eqnarray}
\quad\mathcal{T}_Y(\pi) &=&\sum_i
h_{i}\mathbb{E} \Biggl[ \Biggl| (1-b_1-b_2-b_3)
u^{(0)}_{t_{i }}
\nonumber
\\
&&\hspace*{40pt} {} + \biggl(\frac12-b_1-b_2(1-c_2)-b_3(1-c_3)
\biggr) h_{i}u^{(0,0)}_{t_{i }}
\nonumber
\\
&&\hspace*{40pt} {}+ \biggl(\frac16 - \frac12 b_1 -\frac12
b_2(1-c_2)^2- \frac12b_3(1-c_3)^2
\biggr) h_{i}^2 u^{(0,0,0)}_{t_{i }}
\nonumber
\\[-8pt]
\\[-8pt]
&&\hspace*{40pt} {}- \biggl(b_2 \frac{c_2^2}{2} +
b_3\frac {c_3^2}{2} - b_3 c_2
a_{32} \biggr)h_{i}^2 f^y_{t_{i }}
u^{(0,0)}_{t_{i }}
\nonumber
\\
&&\hspace*{84pt} {} - \biggl(b_2 \frac{c_2^2}{2} +
b_3\frac
{c_3^2}{2} - b_3 c_2 \tilde
\alpha_{32} \biggr)h_{i}^2 \sum
_{\ell=1}^d f^{z^\ell}_{t_{i }}
u^{(\ell,0,0)}_{t_{i }} \Biggr|^2 \Biggr]
\nonumber
\\
&&{} + O \bigl(|\pi|^6 \bigr).
\nonumber
\end{eqnarray}

(2a) If $c_3\neq c_2$, According to steps (1c) and (1d), the
conditions
\begin{eqnarray*}
b_1+b_2+b_3 &=& 1,\qquad
b_2c_2 + b_3c_3 = \tfrac12,
\\
b_2c_2^2 + b_3c_3^2&=&
\tfrac13,\qquad b_3a_{32}c_2 =
b_3 \tilde\alpha_{32}c_2 = \tfrac16
\end{eqnarray*}
%
and
\[
\beta_1 + \beta_2 + \tilde\beta_3 = 1,
\qquad \beta_2c_2 + \tilde \beta_3c_3=
\tfrac12
\]
allow us to obtain an order $3$ method, recalling that $c_2 \neq1$.

%
Observe that the condition on $\beta$ are weaker than on $b$ and that
$a_{32}=\alpha_{32}$. This equality, combined with the other condition
on the coefficients, leads to $a_{jk}=\alpha_{jk}$, $1 \le j,k \le3$.

(2b) Under $(\mathbf Ho)_3$, using the same techniques, as, for
example, in the proof of Theorem~\ref{th s1}, one proves that the above
conditions are necessary.
\section{Four-stage schemes}\label{sec6}

This section is dedicated to the proof of Theorem~\ref{th explicit
barrier}.

We now study the local truncation error for the family of scheme given
by
%
%
\begin{eqnarray}
\label{eq s4 explicit ti2} \qquad Y_{i,2} &=& \mathbb{E}_{t_{i,2}}
\bigl[Y_{i +1} + h_{i}c_{2} f(Z_{i +1})
\bigr],
\\
Z_{i,2} &=& \mathbb{E}_{t_{i,2}} \bigl[H^{\psi_2}_{t_{i,2},c_2h_{i}}
Y_{i +1} + h_{i}c_2 H^{\phi_2}_{t_{i,2},c_2h_{i}}
f( Z_{i +1}) \bigr],
\\
\label{eq s4 explicit ti3} Y_{i,3} &=& \mathbb{E}_{t_{i,3}}
\bigl[Y_{i +1} + h_{i}a_{31} f(Z_{i +1}) +
h_{i} a_{32} f(Z_{i,2}) \bigr],
\\
Z_{i,3} &=& \mathbb{E}_{t_{i,3}} \bigl[H^{\psi_3}_{t_{i,3},c_3h_{i}}
Y_{i +1}
\nonumber
\\[-8pt]
\\[-8pt]
&&\hspace*{20pt} {} + h_{i} \bigl( \alpha_{31}
H^{\phi
_3}_{t_{i,3},c_3h_{i}} f(Z_{i +1}) + \tilde
\alpha_{32} H^{\phi
_3}_{t_{i,3},(c_3-c_2)h_{i}} f(Z_{i,2})
\bigr) \bigr],
\nonumber
\\
\label{eq s4 explicit ti4} Y_{i,4} &=& \mathbb{E}_{t_{i,4}}
\bigl[Y_{i +1} + h_{i}a_{41} f(Z_{i +1}) +
h_{i} a_{42} f(Z_{i,2})+ h_{i}a_{43}
f(Z_{i,3}) \bigr],
\\
Z_{i,4} &=& \mathbb{E}_{t_{i,4}} \bigl[H^{\psi_4}_{t_{i,4},c_4h_{i}}
Y_{i +1}
\nonumber
\\
&&\hspace*{20pt} {} + h_{i} \bigl( \alpha_{41}
H^{\phi_4}_{t_{i,4},c_4h_{i}} f(Z_{i +1}) + \tilde
\alpha_{42} H^{\phi_4}_{t_{i,4},(c_4-c_2)h_{i}} f( Z_{i,2})
\\
&&\hspace*{137pt} {} +\tilde\alpha_{43} H^{\phi_4}_{t_{i,4},(c_4-c_3)h_{i}}
f( Z_{i,3}) \bigr) \bigr].
\nonumber
\end{eqnarray}

The approximation at step~($\mathrm{i}$) is given by
%
%
\begin{eqnarray}
\label{eq RK3 explicit ti67} \quad\qquad Y_{i } &=&
\mathbb{E}_{t_{i }} \bigl[Y_{i
+1} + h_{i} \bigl(
b_{1}f(Z_{i +1})+ b_{2}f(Z_{i,2})+
b_{3}f(Z_{i,3}) + b_{4}f(Z_{i,4})
\bigr) \bigr],
\\
Z_{i } &=& \mathbb{E}_{t_{i }} \bigl[H^{\psi_5}_{t_{i },h_{i}}
Y_{i+1}
\nonumber
\\
&&{}\hspace*{15pt} {} + h_{i} \bigl( \beta_{1}
H^{\phi_5}_{t_{i },h_{i}} f(Z_{i +1})+ \tilde\beta_{2}
H^{\phi_5}_{t_{i },(1-c_2)h_{i}}f( Z_{i,2})
\\
&&\hspace*{39pt} {} + \tilde \beta_{3} H^{\phi_5}_{t_{i },(1-c_3)h_{i}}
f( Z_{i,3}) + \tilde\beta_{4}H^{\phi_5}_{t_{i },(1-c_4)h_{i}}
f( Z_{i,4}) \bigr) \bigr].
\nonumber
\end{eqnarray}

We assume that
\begin{eqnarray*}
a_{31}+a_{32} &=& c_3 \quad\mbox{and}\quad
a_{41}+a_{42}+a_{43} = c_4,
\\
\alpha_{31}+\tilde\alpha_{32} &=& c_3 \quad
\mbox{and}\quad\alpha_{41}+\tilde\alpha_{42}+\tilde
\alpha_{43} = c_4.
\end{eqnarray*}
Moreover, $\psi_2,\psi_3, \psi_4, \psi_5 \in B^3_{[0,1]}$ and
$\phi_2,\phi_3, \phi_4, \phi_5 \in B^2_{[0,1]}$.


We first prove that the following set of condition is necessary to
retrieve an order~$4$ method:

%
\begin{Lemma}\label{le o4 conditions} Assume that $c_2 \neq1$ and
$c_3 \neq1$.
\begin{longlist}[(ii)]
\item[(i)] The order 4 conditions for the $Y$-part are
\begin{eqnarray*}
b_1 + b_2 + b_3 + b_4 &=& 1,
\qquad b_3 \tilde\alpha_{32}c_2 +
b_4 \tilde \alpha_{42} c_2 + b_4
\tilde \alpha_{43} c_3 = \tfrac16,
\\
b_2c_2 + b_3c_3 +
b_4c_4 &=& \tfrac12,\qquad b_3 \tilde
\alpha_{32}c_2c_3 + b_4 \tilde
\alpha_{42} c_2c_3 + b_4 \tilde
\alpha_{43} c_3c_4 = \tfrac18,
\\
b_2c_2^2 + b_3c_3^2
+ b_4c_4^2 &=& \tfrac13,\qquad
b_3 \tilde\alpha_{32}c_2^2 +
b_4 \tilde\alpha_{42} c_2^2 +
b_4 \tilde\alpha_{43} c_3^2 =
\tfrac1{12},
\\
b_2c_2^3 + b_3c_3^3
+ b_4c_4^3 &=& \tfrac14,\qquad
b_4 \tilde\alpha_{43} \tilde\alpha_{32}
c_2 = \tfrac1{24}.
\end{eqnarray*}
%

\item[(ii)] The order 4 conditions for the $Z$-part are
\begin{eqnarray*}
\beta_1 + \beta_2 + \beta_3 &=& 1, \qquad
\beta_2 c_2^2 + \beta_3
c_3^2 = \tfrac13,
\\
\beta_2 c_2 + \beta_3c_3 &=&
\tfrac12,\qquad \beta_3 \alpha_{32} c_3 =
\tfrac16.
\end{eqnarray*}
\end{longlist}
\end{Lemma}

%
%

%
%
\begin{Remark}
(i) If $c_2 = 1$, then $c_3=c_4=1$ and $\beta_1=1$, the approximation
for $Z$ reads
\[
Z_{i } = \mathbb{E}_{t_{i }} \bigl[H^{\psi_5}_{t_{i },h_{i}}
Y_{i +1} + h_{i}H^{\phi_5}_{t_{i },h_{i}}
f(Z_{i +1}) \bigr],
\]
which leads generally to an
order $2$ truncation error for $Z$.

(ii) If $c_2 \neq1$ and $c_3=1$ (then $c_4 = 1$),
\begin{eqnarray*}
Z_{i } &=& \mathbb{E}_{t_{i }} \bigl[H^{\psi_5}_{t_{i },h_{i}}
Y_{i +1} + h_{i}\beta_{1} H^{\phi_5}_{t_{i },h_{i}}
f( Z_{i +1}) +h_{i}\beta_{2} H^{\phi_5}_{t_{i },(1-c_2)h_{i}}
f( Z_{i,2}) \bigr],
\end{eqnarray*}
which leads generally to an order $3$ truncation error for $Z$.
\end{Remark}

\begin{pf*}{Proof of Lemma~\ref{le o4 conditions}}
(1)~We first compute the error expansion at the intermediary steps.
Observe that since we assume that $f$ does not depends on $Y$, we only
need to consider the approximation of $Z$ for the intermediary stages.
\begin{longlist}[(1a)]
\item[(1a)] Error expansion at step 2.

Under $(\mathbf Hr)_4$, using Proposition~\ref{pr weak expansion Z}(i),
we have for $1 \le\ell\le d$,
\begin{eqnarray*}
(\hat{Z}_{t_{i,2}})^\ell&=& \mathbb{E}_{t_{i,2}} \bigl[
\bigl(H^{\psi
_2}_{t_{i,2},c_2h_{i} } \bigr)^\ell Y_{i +1} +
h_{i}c_{2} \bigl(H^{\phi_2}_{t_{i,2},c_2h_{i}}
\bigr)^\ell f({Z}_{t_{i +1}}) \bigr]
\\
&=& u^{(\ell)}_{t_{i,2}} - \frac{c_2^2}{2}h_{i}^2
u^{(\ell,0,0)}_{t_{i,2}} - \frac{c_2^3}{3}h_{i}^3
u^{(\ell,0,0,0)}_{t_{i,2}} + O_{t_{i,2}} \bigl(|\pi|^4
\bigr),
\end{eqnarray*}
which leads to
%
%
\begin{eqnarray}
\label{eq exp bar Z 2} f(\hat{Z}_{t_{i,2}}) &=& -u^{(0)}_{t_{i,2}}
- \frac{c_2^2}{2}h_{i}^2 \sum
_{j=1}^d {}^j v_{t_{i,2}} -
\frac
{c_2^3}{3}h_{i}^3 \sum
_{j=1}^d {}^j w_{t_{i,2}} +
O_{t_{i,2}} \bigl(|\pi|^4 \bigr),
\end{eqnarray}
where we set $ {}^j v = f^{z^j} u^{(j,0,0)}$ and ${}^j w =
f^{z^j} u^{(j,0,0,0)}$, $1 \le j \le d$.


\item[(1b)] Error expansion at step 3.


Observe that, using (\ref{eq exp bar Z 2}), we have for $1 \le\ell
\le d$,
\begin{eqnarray*}
(\hat{Z}_{t_{i,3}})^\ell&=& \mathbb{E}_{t_{i,3}} \bigl[
\bigl(H^{\psi
_3}_{t_{i,3},c_3h_{i} } \bigr)^\ell u_{t_{i +1}}-
h_{i}\alpha_{31} \bigl(H^{\phi _3}_{t_{i,3},c_3h_{i} }
\bigr)^\ell u^{(0)}_{t_{i +1}}
\\
&&\hspace*{74pt} {} - h_{i}\tilde{ \alpha}_{32}
\bigl(H^{\phi_3}_{t_{i,3},(c_3-c_2)h_{i}} \bigr)^\ell u^{(0)}_{t_{i,2}}
\bigr]
\\
&&{}-\mathbb{E}_{t_{i,3}} \Biggl[\tilde{\alpha}_{32}
\frac
{c_2^2}{2}h_{i}^3 \bigl(H^{\phi_3}_{t_{i,3},(c_3-c_2)h_{i}}
\bigr)^\ell\sum_{j=1}^d
{}^j v_{t_{i,2}} \Biggr] +O_{t_{i,3}} \bigl(|
\pi|^4 \bigr).
\end{eqnarray*}
We also used that $\mathbb{E}_{t_{i,3}}
[(H^{\phi_3}_{t_{i,3},(c_3-c_2)h_{i}})^\ell\sum_{j=1}^d {}^j
w_{t_{i,2}} ] = O_{t_{i,3}}(1)$, recalling Proposition~\ref{pr weak
expansion Z} and that under $(\mathbf Hr)_4$, ${}^j w
\in\mathcal{G}^1_b$, $1 \le j \le d$.

Applying Proposition~\ref{pr weak expansion Z}, we compute, recalling
that $\alpha_{31} + \tilde{\alpha}_{32} = c_3$,
%
\begin{eqnarray*}
(\hat{Z}_{t_{i,3}})^\ell&=& u^{(\ell)}_{t_{i,3}} -
\biggl(\frac
{c_3^2}{2}-\tilde{\alpha}_{32}c_2
\biggr)h_{i}^2 u^{(\ell,0,0)}_{t_{i 3}} - \biggl(
\frac{c_3^3}{3}+ \biggl(\frac{c_2^2}{2}-c_2c_3
\biggr)\tilde{\alpha}_{32} \biggr)h_{i} ^3
u^{(\ell,0,0,0)}_{t_{i 3}}
\\
&&{}- \tilde{\alpha}_{32} \frac{c_2^2}{2}h_{i}^3
\mathbb{E}_{t_{i,3}} \Biggl[ \bigl(H^{\phi_3}_{t_{i,3},(c_3-c_2)h_{i}}
\bigr)^\ell\sum_{j=1}^d
{}^j v_{t_{i,2}} \Biggr] + O_{t_{i,3}} \bigl(|
\pi|^4 \bigr).
\end{eqnarray*}

Under $(\mathbf Hr)_4$, ${}^j v \in\mathcal{G}^2_b$, $1 \le j \le d$,
applying Proposition~\ref{pr weak expansion Z}(i), we have that
\[
\mathbb{E}_{t_{i,3}} \bigl[ \bigl(H^{\phi_3}_{t_{i,3},(c_3-c_2)h_{i}}
\bigr)^\ell{}^j v_{t_{i,2}} \bigr] = {}^j
v^{(\ell
)}_{t_{i,3}} + O_{t_{i,2}} \bigl(|\pi| \bigr).
\]

We straightforwardly deduce that
%
%
\begin{eqnarray}
\label{eq exp bar Z 3} f(\hat{Z}_{t_{i,3}}) &=&-u^{(0)}_{t_{i,3}}
- \biggl( \frac{c_3^2}{2}-\tilde{\alpha}_{32}c_2
\biggr)h_{i}^2 \sum_{j=1}^d
{}^j v_{t_{i 3}}
\nonumber
\\
&&{} - \biggl(\frac{c_3^3}{3}+ \biggl( \frac{c_2^2}{2}-c_2c_3
\biggr)\tilde{\alpha }_{32} \biggr)h_{i}^3
\sum_{j=1}^d {}^j
w_{t_{i 3}}
\\
&&{}-\tilde{\alpha}_{32}\frac{c_2^2}{2}h_{i}^3
\sum_{\ell=1}^d \sum
_{j=1}^d f^{z^\ell}_{t_{i,3}}
{}^j v^{(\ell)}_{t_{i,3}} + O_{t_{i,3}} \bigl(|
\pi|^4 \bigr).
\nonumber
\end{eqnarray}
%

\item[(1c)] Error expansion at step 4.

Using (\ref{eq exp bar Z 2})--(\ref{eq exp bar Z 3}), we obtain for $1
\le\ell\le d$,
\begin{eqnarray*}
(\hat{Z}_{t_{i,4}})^\ell&=& \mathbb{E}_{t_{i,4}} \bigl[
\bigl(H^{\psi
_4}_{t_{i,4},c_4h_{i} } \bigr)^\ell u_{t_{i +1}}
\bigr] +O_{t_{i,4}} \bigl(|\pi|^4 \bigr)
\\
&&{} - h_{i}\mathbb{E}_{t_{i,4}} \bigl[\alpha_{41}
\bigl(H^{\phi_4}_{t_{i,4},c_4h_{i}} \bigr)^\ell u^{(0)}_{t_{i +1}}
+ \tilde{\alpha}_{42} \bigl(H^{\phi_4}_{t_{i,4},(c_4-c_2)h_{i}}
\bigr)^\ell u^{(0)}_{t_{i,2}}
\\
&&\hspace*{130pt} {} + \tilde{ \alpha}_{43} \bigl(H^{\phi_4}_{t_{i,4},(c_4-c_3)h_{i}}
\bigr)^\ell u^{(0)}_{t_{i,3}} \bigr]
\\
&&{}- \mathbb{E}_{t_{i,4}} \Biggl[\frac{c_2^2}{2}\tilde{
\alpha}_{42} h_{i}^3 \bigl(H^{\phi_4}_{t_{i,4},(c_4-c_2)h_{i}}
\bigr)^\ell\sum_{j=1}^d
{}^j v_{t_{i,2}}
\\
&&\hspace*{36pt} {} + \biggl(\frac{c_3^2}{2}-\tilde{ \alpha}_{32}c_2
\biggr)\tilde{\alpha}_{43}h_{i}^3
\bigl(H^{\phi_4}_{t_{i,4},(c_4-c_3)h_{i}} \bigr)^\ell\sum
_{j=1}^d {}^j v_{t_{i 3}}
\Biggr].
\end{eqnarray*}
%
%
Using Proposition~\ref{pr weak expansion Z}, recalling that
$\alpha_{41}+\tilde{\alpha}_{42}+\tilde{\alpha}_{43}=c_4$, we compute
\begin{eqnarray*}
(\hat{Z}_{t_{i,4}})^\ell&=&u^{(\ell)}_{t_{i,4}} -
\biggl(\frac{c_4^2}{2} - \tilde{\alpha}_{42}c_2 -
\tilde{\alpha}_{43}c_3 \biggr)h_{i}^2
u^{(\ell,0,0)}_{t_{i 4}}
\\
&&{} - \biggl(\frac{c_4^3}{3}+\tilde{\alpha}_{42}
\frac{c_2^2}{2}- \tilde{\alpha}_{42}c_2c_4
+ \tilde{ \alpha}_{43}\frac{c_3^2}{2}-\tilde{\alpha}_{43}c_3c_4
\biggr)h_{i}^3 u^{(\ell,0,0,0)}_{t_{i 4}}
\\
&&{}+ \tilde{\alpha}_{42} \frac{c_2^2}{2}h_{i}^3
\mathbb{E}_{t_{i,4}} \Biggl[ \bigl(H^{\phi_4}_{t_{i,4},(c_4-c_2)h_{i}}
\bigr)^\ell\sum_{j=1}^d
{}^j v_{t_{i,2}} \Biggr]
\\
&&{}- \tilde{\alpha}_{43} \biggl(\frac{c_3^2}{2}-\tilde{
\alpha}_{32}c_2 \biggr)h_{i} ^3
\mathbb{E}_{t_{i,4}} \Biggl[ \bigl(H^{\phi_4}_{t_{i,4},(c_4-c_3)h_{i}}
\bigr)^\ell\sum_{j=1}^d
{}^j v_{t_{i 3}} \Biggr]
\\
&&{} + O_{t_{i,4}} \bigl(|\pi|^4 \bigr).
\end{eqnarray*}

Applying Proposition~\ref{pr weak expansion Z}, this leads to,
recalling that $(\mathbf Hr)_4$ is in force,
%
%
\begin{eqnarray}
\label{eq exp bar Z 4} \qquad f(\hat{Z}_{t_{i,4}}) &=& -u^{(0)}_{t_{i,4}}-
\biggl( \frac{c_4^2}{2} - \tilde{\alpha}_{42}c_2 -
\tilde{ \alpha}_{43}c_3 \biggr)h_{i}^2
\sum_{j=1}^d {}^j
v_{t_{i 4}}
\nonumber
\\
&&{} - \biggl(\frac{c_4^3}{3}+\tilde{\alpha}_{42}
\frac{c_2^2}{2}- \tilde{\alpha}_{42}c_2c_4
+ \tilde{ \alpha}_{43}\frac{c_3^2}{2}-\tilde{\alpha}_{43}c_3c_4
\biggr)h_{i}^3 \sum_{j=1}^d
{}^j w_{t_{i 4}}
\nonumber
\\[-8pt]
\\[-8pt]
&&{} - \biggl( \tilde{\alpha}_{42} \frac{c_2^2}{2} +\tilde{
\alpha}_{43} \biggl(\frac{c_3^2}{2}-\tilde{\alpha}_{32}c_2
\biggr) \biggr) h_{i}^3 \sum
_{\ell=1}^d \sum_{j=1}^d
f^{z^\ell}_{t_{i,4}} {}^j v^{(\ell)}_{t_{i,4}}
\nonumber
\\
&&{}+ O_{t_{i,4}} \bigl(| \pi|^4 \bigr).
\nonumber
\end{eqnarray}

\item[(2a)] We now study the error for the $Y$-part at the
final step.

Using (\ref{eq exp bar Z 2})--(\ref{eq exp bar Z 4}), we obtain
\begin{eqnarray*}
\hat{Y}_{t_{i }} &=& \mathbb{E}_{t_{i }} \bigl[u_{t_{i +1}} -
h_{i} \bigl( b_1 u^{(0)}_{t_{i +1}} +
b_2 u^{(0)}_{t_{i,2}} + b_3
u^{(0)}_{t_{i,3}} + b_4 u^{(0)}_{t_{i,4}}
\bigr) \bigr]
\\
&&{}- b_2\mathbb{E}_{t_{i }} \Biggl[\frac{c_2^2}{2}h_{i}^3
\sum_{j=1}^d {}^j
v_{t_{i,2}} + \frac{c_2^3}{3}h_{i}^4 \sum
_{j=1}^d {}^j
w_{t_{i,2}} \Biggr]
\\
&&{} -b_3\mathbb{E}_{t_{i }} \Biggl[ \biggl(
\frac{c_3^2}{2}- \tilde{\alpha}_{32}c_2
\biggr)h_{i}^3 \sum_{j=1}^d
{}^j v_{t_{i 3}} + \biggl(\frac{c_3^3}{3}+\tilde{
\alpha}_{32} \frac{c_2^2}{2}-\tilde{\alpha}_{32}c_2c_3
\biggr)h_{i}^4 \sum_{j=1}^d
{}^j w_{t_{i 3}}
\\
&&\hspace*{207pt} {} + \tilde{\alpha}_{32}\frac{c_2^2}{2}h_{i}^4
\sum_{\ell=1}^d \sum
_{j=1}^d f^{z^\ell}_{t_{i,3}}
{}^j v^{(\ell)}_{t_{i,3}} \Biggr]
\\
&&{} -b_4\mathbb{E}_{t_{i }} \Biggl[ \biggl(
\frac{c_4^2}{2} - \tilde{\alpha}_{42}c_2 - \tilde{\alpha
}_{43}c_3 \biggr)h_{i}^3 \sum
_{j=1}^d {}^j
v_{t_{i 4}}
\\
&&\hspace*{41pt} {} + \biggl(\frac
{c_4^3}{3}+\tilde{\alpha}_{42}
\frac{c_2^2}{2}-\tilde{\alpha}_{42}c_2c_4 +
\tilde{\alpha}_{43}\frac{c_3^2}{2}-\tilde{\alpha
}_{43}c_3c_4 \biggr)h_{i}^4
\sum_{j=1}^d {}^j
w_{t_{i 4}} \Biggr]
\\
&&{} -b_4\mathbb{E}_{t_{i }} \Biggl[ \biggl( \tilde{
\alpha}_{42} \frac{c_2^2}{2} + \tilde{\alpha}_{43}
\biggl( \frac{c_3^2}{2}-\tilde{\alpha}_{32}c_2 \biggr)
\biggr) h_{i}^4 \sum_{\ell=1}^d
\sum_{j=1}^d f^{z^\ell
}_{t_{i,4}}
{}^j v^{(\ell)}_{t_{i,4}} \Biggr]
\\
&&{} +O_{t_{i }} \bigl(|\pi|^5 \bigr).
\end{eqnarray*}

Under $(\mathbf Hr)_4$,
using Proposition~\ref{pr weak expansion Z}, we then compute 
%
\begin{eqnarray*}
\hat{Y}_{t_{i }} &=& u_{t_{i }} + h_{i}(1-b_1-b_2-b_3-b_4
)u^{(0)}_{t_{i }}
\\
&&{} +h_{i}^2 \biggl( \frac{1}{2} -
b_2(1-c_2)- b_3(1-c_3)-
b_4(1-c_4) \biggr) u^{(0,0)}_{t_{i }}
\\
&&{} + h_{i}^3 \biggl( \frac{1}{6}-
\frac
{b_2(1-c_2)^2+b_3(1-c_3)^2+b_4(1-c_4)^2}{2} \biggr) u^{(0,0,0)}_{t_{i }}
\\
&&{} + h_{i}^4 \biggl( \frac{1}{24}-
\frac
{b_2(1-c_2)^3+b_3(1-c_3)^3+b_4(1-c_4)^3}{6} \biggr) u^{(0,0,0,0)}_{t_{i }}
\\
&&{} - h_{i}^3 \biggl( b_2\frac{c_2^2}{2}
+ b_3 \biggl(\frac{c_3^2}{2}-\tilde{\alpha}_{32}c_2
\biggr) + b_4 \biggl(\frac{c_4^2}{2} - \tilde{
\alpha}_{42}c_2 - \tilde{\alpha}_{43}c_3
\biggr) \biggr)\sum_{j=1}^d
{}^j v_{t_{i }}
\\
&&{} - h_{i}^4 \biggl( b_2
\frac{c_2^2}{2}(1-c_2) + b_3 \biggl(
\frac
{c_3^2}{2}-\tilde{\alpha}_{32}c_2 \biggr)
(1-c_3)
\\
&&\hspace*{51pt} {}+ b_4 \biggl(\frac{c_4^2}{2} - \tilde{
\alpha}_{42}c_2 - \tilde{\alpha}_{43}c_3
\biggr) (1-c_4) \biggr)\sum_{j=1}^d
{}^j v^{(0)}_{t_{i }}
\\
&&{} - h_{i}^4 \biggl( b_2\frac{c_2^3}{3}
+ b_3 \biggl(\frac{c_3^3}{3}+\tilde{\alpha}_{32}
\frac{c_2^2}{2}-\tilde{\alpha}_{32}c_2c_3
\biggr)
\\
&&\hspace*{28pt} {}+ b_4 \biggl(\frac{c_4^3}{3}+\tilde{
\alpha}_{42} \frac{c_2^2}{2}-\tilde{\alpha}_{42}c_2c_4
+ \tilde{\alpha}_{43}\frac{c_3^2}{2}-\tilde{\alpha}_{43}c_3c_4
\biggr) \biggr) \sum_{j=1}^d
{}^j w_{t_{i }}
\\
&&{} - h_{i}^4 \biggl( b_3\tilde{
\alpha}_{32} \frac{c_2^2}{2}+ b_4 \biggl( \tilde{
\alpha}_{42} \frac{c_2^2}{2} + \tilde{\alpha}_{43}
\biggl( \frac{c_3^2}{2}-\tilde{\alpha}_{32}c_2 \biggr)
\biggr) \biggr) \sum_{\ell=1}^d\sum
_{j=1}^d f^{z^\ell}_{t_{i }}
{}^j v^{(\ell)}_{t_{i }}
\\
&&{} +O_{t_{i }} \bigl(|\pi|^5 \bigr).
\end{eqnarray*}

Under $(\mathbf Ho)_4$, using the same techniques as in the proof of
Theorem~\ref{th 2 stage scheme imp}, one proves inductively on the
order that each factor has to be equal to $0$, which leads to the set
(i) of conditions of the lemma. It appears that these conditions are
the same as in the ODE case. From, for example, Section~322, page 175
in \cite{but08}, we know that $c_4=1$ necessarily.

\item[(3)] We now study the error for the $Z$-part at the final step, taking
into account $c_4 = 1$ and $c_2<c_3<1$. 
We thus have
\begin{eqnarray*}
\hat{Z}_{t_{i }} &=& \mathbb{E}_{t_{i }} \bigl[H^{\psi_5}_{t_{i },h_{i}}
{Y}_{t_{i +1}} + h_{i} \bigl( \beta_{1}
H^{\phi_5}_{t_{i },h_{i}} f({Z}_{t_{i +1}})+ \beta_{2}H^{\phi
_5}_{t_{i },(1-c_2)h_{i}}f(
\hat{Z}_{t_{i,2}})
\\
&&\hspace*{162pt} {} + \tilde\beta_{3}H^{\phi
_5}_{t_{i },(1-c_3)h_{i}}
f({Z}_{t_{i,3}}) \bigr) \bigr].
\end{eqnarray*}
We are thus considering a 3-stage scheme for the $Z$ part. Using the
results of step~1, we obtain the following
expansion, for $1 \le\ell\le d$: 
\begin{eqnarray*}
\hat{Z}_{t_{i }}^\ell&=& {Z}_{t_{i }}^\ell+
(1-\beta_1-\beta_2-\beta_3)
u^{(\ell,0)}_{t_{i }}
\\
&&{} + \biggl(\frac12-\beta_1- \beta_2(1-c_2)-
\beta_3(1-c_3) \biggr) h_{i}u^{(\ell,0,0)}_{t_{i }}
\\
&&{} + \biggl(\frac16 -\frac12 \beta_1 -\frac12
\beta_2(1-c_2)^2- \frac12
\beta_3(1-c_3)^2 \biggr)
h_{i}^2 u^{(\ell,0,0,0)}_{t_{i }}
\\
&&{} - \biggl(\beta_2 \frac{c_2^2}{2} + \beta_3
\frac{c_3^2}{2} - \beta_3 c_2 \tilde
\alpha_{32} \biggr)h_{i}^2 f^{z^{\ell}}_{t_{i }}
\sum_{j=1}^d u^{(j,0,0)}_{t_{i }}
+ O_{t_{i }} \bigl(|\pi|^3 \bigr).
\end{eqnarray*}

It is then obvious that set (ii) of the condition is sufficient to
obtain an order $4$ truncation error on $Z$. Moreover, arguing as, for
example, in steps (2b)--(2c) of the proof of Theorem~\ref{th 2 stage
scheme imp}, by induction on the order required, one proves that these
condition are also necessary, provided that $(\mathbf Ho)_4$ is in
force.
\end{longlist}\upqed
\end{pf*}


%
%

\begin{pf*}{Proof of Theorem~\ref{th explicit barrier}}
The set of condition (ii) leads, using case I of Theorem~\ref{co 3 stage}, with
$(b_j)=(\beta_j)$ and $(a_{kj}) = (\alpha_{kj})$, to the only possible
value for $\alpha_{32}$ is given by
\[
\alpha_{32} = \frac{c_3(c_3-c_2)}{c_2(2-3c_2)}.
\]
In our context equations (322b) and (322c) in \cite{but08} read
\begin{eqnarray*}
b_4\tilde\alpha_{43} (c_3-c_2)
c_3 &=& \frac1{12} - \frac{c_2}{6},
\\
b_4\tilde\alpha_{43} \tilde\alpha_{32}
c_2 &=& \tfrac1{24}.
\end{eqnarray*}
Dividing these two equations, we obtain
\[
\frac{(c_3-c_2) c_3}{\alpha_{32} c_2 } = 2 - 4 c_2.
\]
It follows from the expression of $\alpha_{32}$ that $c_2=0$, which is
not possible.
\end{pf*}


\setcounter{equation}{0}
%
\begin{appendix}\label{sec7}
\section*{Appendix}

\subsection{Schemes stability}

\subsubsection{\texorpdfstring{Proof of Theorem \protect\ref{th stab}}
{Proof of Theorem 1.1}}\label{sec7.1.1}
Using (\ref{eq prop PhiY})--(\ref{eq prop PhiZ}), we compute, for $1
\ge\eta>0$ to be fixed later on that
%
%
\begin{eqnarray}
\label{eq stab explicit 1} |\delta Y_{i }|^{2} &\le&
\biggl(1+\frac{h_{i}}{\eta} \biggr) \bigl|\mathbb{E}_{t_{i }} [\delta
Y_{i +1} ] \bigr|^{2} + C \frac{\eta}{h_{i}} \bigl|\mathbb
{E}_{t_{i }} \bigl[\zeta^{Y}_{i } \bigr]
\bigr|^{2}
\nonumber
\\
&&{}+ Ch_{i}^{2} \biggl(1+ \frac{\eta}{h_{i}} \biggr)
\biggl(\frac1h B_{i } + \mathbb{E}_{t_{i }} \bigl[|\delta
Y_{i +1}|^2+|\delta Z_{i
+1}|^2 \bigr]
\biggr),
\\
|\delta Z_{i }|^{2} &\le& C \biggl(\frac1{h_{i}}
B_{i } +h_{i}\mathbb{E}_{t_{i }} \bigl[|\delta
Y_{i +1}|^2 + |\delta Z_{i +1}|^2
\bigr] + \bigl|\mathbb{E}_{t_{i }} \bigl[\zeta^{Z}_{i }
\bigr] \bigr|^{2} \biggr),
\nonumber
\end{eqnarray}
where $B_{i }:= \mathbb{E}_{t_{i }} [|\delta Y_{i +1}|^2 -
|\mathbb{E}_{t_{i }} [\delta Y_{i +1} ]|^{2} ]$.\vadjust{\goodbreak}

Defining for $1 \ge\varepsilon> 0$ to be fixed later on $I_{i
}^{\varepsilon}:= |\delta Y_{i }|^{2} + \varepsilon h_{i}|\delta Z_{i
}|^{2} $, we compute
\begin{eqnarray*}
&& I_{i }^{\varepsilon/2} + \frac{\varepsilon}{2} h_{i}|
\delta Z_{i }|^2
\\
&&\qquad \le \biggl(1+\frac{h_{i}}{\eta} \biggr) \bigl| \mathbb{E}_{t_{i }} [
\delta Y_{i
+1} ] \bigr|^{2} + C( \varepsilon+\eta)
B_{i } + C\frac{\eta}{h_{i}} \bigl|\mathbb{E}_{t_{i }} \bigl[
\zeta^{Y}_{i } \bigr] \bigr|^{2} + C h_{i}
\bigl|\mathbb{E}_{t_{i }} \bigl[\zeta^{Z}_{i } \bigr]
\bigr|^{2}
\\
&&\quad\qquad{} + \biggl(Ch_{i}^{2} \biggl(1+
\frac{\eta}{h_{i}} \biggr) + C \varepsilon h_{i}^{2} \biggr)
\mathbb{E}_{t_{i }} \bigl[|\delta Y_{i,+1}|^{2} + |
\delta Z_{i,+1}|^{2} \bigr]. 
\end{eqnarray*}
Setting $\varepsilon=\eta=\frac1{2C}$ and observing that
$|\mathbb{E}_{t_{i }} [\delta Y_{i +1} ]|^{2} = \mathbb{E}_{t_{i }}
[|\delta Y_{i +1}|^2 ] - B_{i }$, we compute that, for $h^*$ small
enough
%
%
\begin{eqnarray}
\label{eq pre final} \qquad\quad I_{i }^{\varepsilon/2} +
\frac{\varepsilon}{2} h_{i}|\delta Z_{i }|^2 &\le&
(1+Ch_{i})I_{i +1}^{\varepsilon/2} + C\frac{\eta}{h_{i}} \bigl|
\mathbb{E}_{t_{i }} \bigl[\zeta^{Y}_{i } \bigr]
\bigr|^{2} + C h_{i} \bigl|\mathbb{E}_{t_{i }} \bigl[
\zeta^{Z}_{i } \bigr] \bigr|^{2}.
\end{eqnarray}
Using the discrete version of Gronwall's lemma, we obtain
\begin{eqnarray*}
\max_{0\le i \le n- 1} \mathbb{E} \bigl[| \delta Y_{i }|^{2}
\bigr] &\le& \max_{0\le i \le n-1} I_{i }^{\varepsilon/2}
\\
&\le& C \Biggl( I_n^{\varepsilon/2} + \sum
_{i=0}^{n-1} h_{i}\mathbb{E} \biggl[
\frac{1}{h_{i}^{2}} \bigl|\mathbb{E}_{t_{i }} \bigl[\zeta^{Y}_{i }
\bigr] \bigr|^{2} + \bigl| \mathbb{E}_{t_{i }} \bigl[
\zeta^{Z}_{i } \bigr] \bigr|^{2} \biggr] \Biggr).
\end{eqnarray*}
The control of $\sum_{i=0}^{n-1}h_i \mathbb{E} [|\delta Z_{i
}|^{2} ]$ is then
obtained summing inequality (\ref{eq pre final}) over~$i$. 

\subsubsection{\texorpdfstring{Proof of Proposition \protect\ref{pr pre-conv}}
{Proof of Proposition 1.1}}\label{sec7.1.2}

We simply observe that the solution $(Y,Z)$ of the BSDE is also the
solution of a perturbed scheme with $\zeta^{Y}_{i }:=
\hat{Y}_{t_{i }}-Y_{t_{i }}$ and $\zeta^{Z}_{i }:= \hat{Z}_{t_{i
}}-Z_{t_{i }}$,
and with terminal conditions $\tilde{Y}_n:= g(X_T)$ and
$\tilde{Z}_n:=\nabla g^\top(X_T)\sigma(X_T)$. The proof then follows
directly from Theorem~\ref{th stab}. 

\subsubsection{\texorpdfstring{Proof of Theorem \protect\ref{th conv RK scheme}}
{Proof of Theorem 1.2}}\label{sec7.1.3}
Claim (ii) is a direct application of (i) and
Proposition~\ref{pr ge conv res}.

We now prove (i).
\begin{longlist}[(2a)]
\item[(1)] We define $U_{i,j}$ (resp., $\tilde{U}_{i,j}$) and $V_{i,j}$
(resp., $\tilde{V}_{i,j}$) as $Y_{i,j}$ and $Z_{i,j}$ in
Definition~\ref{de RK scheme}(ii) using $U$ (resp., $\tilde{U}$)
instead of $Y_{i +1}$ and $V$ (resp., $\tilde{V}$) instead of $Z_{i
+1}$. Let us also denote
\[
F_{i,j}:= f(U_{i,j},V_{i,j}),\qquad
\tilde{F}_{i,j}:= f(\tilde{U}_{i,j},\tilde{V}_{i,j}) \quad\mbox{and}\quad\delta F_{i,j}:= F_{i,j}-
\tilde{F}_{i,j}.
\]
With this notation, we have that
\[
\Phi_i^Y(U,V):= \sum_{j=1}^{q+1}
b_{j}f(U_{i,j},V_{i,j}) \quad\mbox{and}\quad
\Phi_i^Y( \tilde U, \tilde V) = \sum
_{j=1}^{q+1} b_{j}f(\tilde{U}_{i,j},
\tilde{V}_{i,j}).
\]
Since $f$ is Lipschitz-continuous, we compute
\begin{eqnarray*}
\mathbb{E}_{t_{i }} \bigl[ \bigl| \Phi_i^Y(U,V) -
\Phi_i^Y(\tilde U, \tilde V) \bigr|^2 \bigr] & \le&
C \mathbb{E}_{t_{i }} \bigl[|\delta{U}_{i,1}|^2+|
\delta{V}_{i,1}|^2 \bigr] + \sum
_{j=2}^{q+1} \mathbb{E}_{t_{i }} \bigl[|\delta
F_{i,j}|^2 \bigr].
\end{eqnarray*}
We also have that
\[
\mathbb{E}_{t_{i }} \bigl[\Phi_i^Z(U,V) \bigr] =
\sum_{j=1}^{q} \beta_{j}
H^i_{q+1,j} f(U_{i,j},V_{i,j})
\]
and
\[
\mathbb{E}_{t_{i }} \bigl[\Phi_i^Z( \tilde U,
\tilde V) \bigr] = \sum_{j=1}^{q}
\beta_{j} H^i_{q+1,j} f(\tilde{U}_{i,j},
\tilde{V}_{i,j}).
\]
Combining the Cauchy--Schwarz inequality with property (\ref{eq prop
H}) and the Lipschitz continuity of $f$, we compute
%
\begin{eqnarray*}
&& \mathbb{E}_{t_{i }}
\bigl[ \bigl| \Phi_i^Z(U,V) - \Phi_i^Z(
\tilde U, \tilde V) \bigr|^2 \bigr]
\\
&&\qquad \le \frac{C}{h_{i}} \mathbb{E}_{t_{i }} \bigl[|
\delta{U}_{i,1}|^2+| \delta{V}_{i,1}|^2
\bigr] +\frac{C}{h_{i}} \sum_{j=2}^q
\mathbb{E}_{t_{i }} \bigl[|\delta F_{i,j}|^2 \bigr].
\end{eqnarray*}
Moreover, we observe, using the Lipschitz-continuity property of $f$,
\[
\mathbb{E}_{t_{i }} \bigl[|\delta F_{i,j}|^2 \bigr]
\le C\mathbb{E}_{t_{i }} \bigl[|\delta{U}_{i,j}|^2 +
|\delta{V}_{i,j}|^2 \bigr].
\]
%

\item[(2a)] For $j=2$, we compute that
\begin{eqnarray*}
\mathbb{E}_{t_{i }} \bigl[|\delta{U}_{i,2}|^2
\bigr] &\le& C \bigl( \mathbb{E}_{t_{i }} \bigl[|\delta
{U}_{i,1}|^2+h_{i}^2 |
\delta{V}_{i,1}|^2 \bigr] + h_{i}^2
\mathbb{E}_{t_{i }} \bigl[|\delta{U}_{i,2}|^2+|
\delta{V}_{i,2}|^2 \bigr] \bigr),
\\
\mathbb{E}_{t_{i }} \bigl[|\delta{V}_{i,2}|^2
\bigr] &\le& C \biggl(\frac1h_{i}\mathbb{E}_{t_{i }} \bigl[|
\delta{U}_{i,1}|^2 - \bigl|\mathbb{E}_{t_{i }} [
\delta{U}_{i,1} ] \bigr|^2 \bigr]+ h_{i}\mathbb
{E}_{t_{i }} \bigl[|\delta{U}_{i,1}|^2+|
\delta{V}_{i,1}|^2 \bigr] \biggr).
\end{eqnarray*}

For $|\pi|$ small enough, we then obtain
\begin{eqnarray*}
&& \mathbb{E}_{t_{i }} \bigl[|\delta{U}_{i,2}|^2+|
\delta{V}_{i,2}|^2 \bigr]
\\
&&\qquad \le C \biggl( \frac1h_{i}\mathbb{E}_{t_{i }} \bigl[|
\delta{U}_{i,1}|^2 - \bigl|\mathbb{E}_{t_{i }} [
\delta{U}_{i,1} ] \bigr|^2 \bigr]+ \mathbb{E}_{t_{i }}
\bigl[|\delta{U}_{i,1}|^2+ h_{i}|
\delta{V}_{i,1}|^2 \bigr] \biggr),
\end{eqnarray*}
which, since $f$ is Lipschitz, straightforwardly leads to
\[
\mathbb{E}_{t_{i }} \bigl[|\delta F_{i,2}|^2 \bigr]
\le C \biggl(\frac1h_{i}\mathbb{E}_{t_{i }} \bigl[|
\delta{U}_{i,1}|^2 - \bigl|\mathbb{E}_{t_{i }} [
\delta{U}_{i,1} ] \bigr|^2 \bigr]+ \mathbb{E}_{t_{i }}
\bigl[|\delta{U}_{i,1}|^2+ h_{i}|
\delta{V}_{i,1}|^2 \bigr] \biggr).
\]

\item[(2b)] For $2<j\le q+1$, we have that
\begin{eqnarray*}
\mathbb{E}_{t_{i }} \bigl[|\delta{U}_{i,j}|^2
\bigr] &\le& C\mathbb{E}_{t_{i }} \Biggl[|\delta{U}_{i,1}|^2+h_{i}^2
\sum_{k=1}^{j-1} |\delta
F_{i,j}|^2 \Biggr] + C h_{i}^2
\mathbb{E}_{t_{i }} \bigl[|\delta{U}_{i,j}|^2+|
\delta{V}_{i,j}|^2 \bigr],
\\
\mathbb{E}_{t_{i }} \bigl[|\delta{V}_{i,j}|^2
\bigr] &\le& C \Biggl(\frac1h_{i}\mathbb{E}_{t_{i }} \bigl[|
\delta{U}_{i,1}|^2 - \bigl|\mathbb{E}_{t_{i }} [
\delta{U}_{i,1} ] \bigr|^2 \bigr]+ h_{i}\sum
_{k=1}^{j-1}\mathbb{E}_{t_{i }} \bigl[|\delta
F_{i,j}|^2 \bigr] \Biggr)
\end{eqnarray*}
and for $|\pi|$ small enough,
\begin{eqnarray*}
&& \mathbb{E}_{t_{i }} \bigl[|\delta F_{i,j}|^2
\bigr]
\\
&&\qquad \le C \Biggl(\frac1h_{i} \mathbb{E}_{t_{i }} \bigl[|
\delta{U}_{i,1}|^2 - \bigl| \mathbb{E}_{t_{i }} [
\delta{U}_{i,1} ] \bigr|^2 \bigr]+ \mathbb{E}_{t_{i }}
\bigl[|\delta{U}_{i,1}|^2 \bigr]+h_{i}\sum
_{k=1}^{j-1} \mathbb{E}_{t_{i }}
\bigl[|\delta F_{i,j}|^2 \bigr] \Biggr).
\end{eqnarray*}
An easy mathematical induction completes the proof. 
\end{longlist}

\subsection{It\^o--Taylor expansions}\label{sec7.2}

\subsubsection{\texorpdfstring{Proof of Proposition \protect\ref{pr weak expansion Y}}
{Proof of Proposition 2.2}}\label{sec7.2.1}
Using Proposition~\ref{pr basic expansion} (Theorem
5.5.1 in \cite{klopla92}), we compute
\[
v(t+h,X_{t+h}) = \sum_{\alpha\in\mathcal{A}_m}v^\alpha_tI^{\alpha}_{t,t+h}
+ \sum_{\beta\in\mathcal{A}_{m+1} \setminus\mathcal{A}_m}I^{\beta
}_{t,t+h}
\bigl[v^\beta \bigr],
\]
recalling that $\mathcal{B}(\mathcal{A}_{m}) = \mathcal{A}_{m+1}
\setminus\mathcal{A}_{m} $.

Taking the conditional expectation on both sides and using Lemma 5.7.1
in \cite{klopla92},
we obtain
\[
\Biggl\llvert\mathbb{E}_{t} \bigl[v(t+h,X_{t+h}) \bigr]-
\sum_{k=0}^m v^{(0)_k}_t
\frac
{h^k}{k!} \Biggr\rrvert= \bigl\llvert\mathbb{E}_{t}
\bigl[I^{(0)_{m+1}}_{t,t+h} \bigl[v_{\cdot}^{(0)_{m+1}}
\bigr] \bigr] \bigr\rrvert.
\]
Since $v \in\mathcal{G}_b^{\beta}$ for all $\beta\in\mathcal
{A}_{m+1}$, in
particular $v \in\mathcal{G}_b^{(0)_{m+1}}$, we obtain
\[
\bigl\llvert\mathbb{E}_{t} \bigl[I^{(0)_{m+1}}_{t,t+h}
\bigl[v^{(0)_{m+1}} \bigr] \bigr] \bigr\rrvert= O_t
\bigl(h^{m+1} \bigr),
\]
which completes the proof.

\subsubsection{\texorpdfstring{Proof of Proposition \protect\ref{pr weak expansion Z}}
{Proof of Proposition 2.3}}\label{sec7.2.2}
(i) (1) Using Proposition~\ref{pr basic expansion} (Theorem 5.5.1 in
\cite{klopla92}), we compute
\begin{eqnarray*}
&& \bigl(H^\psi_{t,t+h} \bigr)^\ell
v(t+h,X_{t+h}) - \sum_{\alpha\in\mathcal{A}
_{m+1}}v^\alpha_t
\bigl(H^\psi_{t,t+h} \bigr)^\ell
I^{\alpha}_{t,t+h}
\\
&&\qquad = \sum_{\beta\in\mathcal{A}_{m+2} \setminus
\mathcal{A}_{m+1}} \bigl(H^\psi_{t,t+h}
\bigr)^\ell I^{\beta}_{t,t+h} \bigl[v_{\cdot}^\beta
\bigr],
\end{eqnarray*}
recalling that $\mathcal{B}(\mathcal{A}_{m+1}) = \mathcal{A}_{m+2}
\setminus\mathcal{A}_{m+1} $.\vadjust{\goodbreak}

(2) We now compute $\mathbb{E}_{t} [(H^\psi_{t,t+h})^\ell
I^{\alpha}_{t,t+h} ]$ for $\alpha\in\mathcal{A}_{m+1}$, recalling
that\break  $(H^\psi_{t,t+h})^\ell:= \frac1h
I^{(\ell)}_{t,t+h}[\psi^\ell_{t,h}]$; see Definition~\ref{de
psi-H}(ii).

If $\alpha^+ \neq(\ell)$, we observe that $\mathbb{E}_{t}
[(H^\psi_{t,t+h})^\ell I^{\alpha}_{t,t+h} ]=0$; see, for example,
Lemma~5.7.2 in \cite{klopla92}.

Now, let $\alpha$ be such that $\ell(\alpha)=q$, $1\le q\le m+1$ and
$\alpha^+ = (\ell)$. Then there exits $1\le l \le q$, such that
$\alpha= (0)_{l-1}*(\ell)*(0)_{q-l}$, and we have
%
%
\begin{eqnarray}
\label{eq weak z 1} \qquad && \mathbb{E}_{t} \bigl[
\bigl(H^\psi_{t,h} \bigr)^\ell
I^{\alpha}_{t,t+h} \bigr]
\nonumber
\\
&&\qquad = \frac1h I^{(0)_{q-l}}_{t,t+h} \biggl[
\mathbb{E}_{t} \biggl[I^{(\ell
)}_{t,\cdot} \biggl[
\psi^\ell \biggl( \frac{\cdot-t}{h} \biggr) \biggr]I^{(\ell
)}_{t,\cdot}
\bigl[I^{(0)_{l-1}}_{t,\cdot} \bigr] \biggr] \biggr]
\nonumber
\\
&&\qquad = \frac1{h(l-1)!}I^{(0)_{q-l}}_{t,t+h} \biggl[
I^{(0)}_{t,\cdot} \biggl[\psi^\ell \biggl(
\frac{\cdot-t}{h} \biggr) (\cdot-t)^{l-1} \biggr] \biggr]
\\
&&\qquad = \frac1{h(l-1)!(q-l)!}\int_t^{t+h}
(t+h-u)^{q-l}(u-t)^{l-1} \psi^\ell \biggl(
\frac{u-t}{h} \biggr) \,\mathrm{d}u
\nonumber
\\
&&\qquad =\frac{h^{q-1}}{(l-1)!(q-l)!}\int_0^{1}
(1-r)^{q-l}r^{l-1} \psi^\ell(r) \,\mathrm{d}r.
\nonumber
\end{eqnarray}

Since $\psi^\ell\in\mathcal{B}^m_{[0,1]}$,
\[
\mathbb{E}_{t} \bigl[ \bigl(H^\psi_{t,h}
\bigr)^\ell I^{\alpha}_{t,t+h} \bigr] =
\frac
{h^{q-1}}{(q-1)!} \mathbf{1}_{\{\alpha_1=\ell\}}.
\]

(3) Using Lemma 5.7.2 in \cite{klopla92}, for $\beta\in
\mathcal{A}_{m+2} \setminus\mathcal{A}_{m+1} $ and $1\le j \le d$, we
have
\[
\mathbb{E}_{t} \bigl[ \bigl(H^\psi_{t,t+h}
\bigr)^j I^{\beta}_{t,t+h} \bigl[v^\beta
\bigr] \bigr]:=\frac1h\mathbb{E}_{t} \bigl[I^{(j)}_{t,t+h}
\bigl[\psi^j_{t,h} \bigr]I^{\beta
}_{t,t+h}
\bigl[v^\beta \bigr] \bigr] = 0\qquad\mbox{if } \beta^+ \neq(j).
\]
We are now considering $\beta\in\mathcal{A}_{m+2} \setminus\mathcal
{A}_{m+1} $ such that $\beta^+ = (j)$, that is, $\beta$ with at most
one nonzero component. According to the notation of Lemma 5.7.2 in
\cite{klopla92} (see the beginning of Section~5.7 in \cite{klopla92}),
we then compute that
\[
k_0(\beta) + k_1(\beta) = m+1 \quad\mbox{and}\quad
k_0 \bigl((j) \bigr) = k_1 \bigl((j) \bigr) = 0.
\]
Since $\ell((j)^+)=1$, we obtain $\omega((j),\beta) = m+2$ and using
again Lemma 5.7.2, we obtain
\[
\biggl|\mathbb{E}_{t} \biggl[\sum_{\beta\in\mathcal{A}_{m+2}}
\bigl(H^\psi_{t,t+h} \bigr)^j I^{\beta}_{t,t+h}
\bigl[v^\beta \bigr] \biggr] \biggr| = O_t
\bigl(h^{m+1} \bigr),
\]
recalling that $v \in\mathcal{G}_b^\beta$, for $\beta\in\mathcal
{A}_{m+2}$.

(ii) This is a straightforward consequence of It\^o's
formula applied to $v$ and the fact that $v^{(0)}$ and $v^{(\ell)}$ are
bounded under $\mathcal{G}^1_b$.

(iii) We follow the arguments of (i). In particular, since
$\psi=(1,\ldots,1)$ in (\ref{eq weak z 1}), using the basic properties
of the Beta function, one obtains
\[
\mathbb{E} \bigl[ \bigl(H^\psi_{t,h} \bigr)^j
I^{\alpha}_{t,t+h} \bigr] = \frac
{h^{q-1}}{q!}
\]
for $\ell(\alpha)=q$ and $\alpha^+ = (j)$, $1\le q\le m+1$, $ 1\le j
\le d$. The proof is completed observing that $v^{(\alpha)} =
v^{(j)*(0)_{q-1}}$ for such $\alpha$ under the assumption
$L^{(0)}\circ
L^{(j)} = L^{(j)}\circ L^{(0)}$.
\end{appendix}




\printaddresses

\end{document}